\documentclass[12pt]{amsart}
\usepackage{amsmath,amscd,amssymb,amsfonts,graphics}
\usepackage{xcolor,hyperref}
\hypersetup{colorlinks=true,citecolor=black,linkcolor=black}
\newcommand{\hl}{\hyperlink}
\newcommand{\htt}{\hypertarget}
\newcommand{\h}{\hbox}
\newcommand{\q}{\quad}
\newcommand{\nin}{\noindent}
\newcommand{\bsn}{\par\bigskip\noindent}
\newcommand{\msn}{\par\medskip\noindent}
\newcommand{\skn}{\par\smallskip\noindent}
\newcommand{\bs}{\par\bigskip}
\newcommand{\ms}{\par\medskip}
\newcommand{\sk}{\par\smallskip}
\newcommand{\mprod}{\h{$\prod$}}
\newcommand{\sotim}{\h{$\otimes$}}

\newcommand{\msum}{\h{$\sum$}}
\newcommand{\mopl}{\h{$\bigoplus$}}
\newcommand{\ssb}{\raise.15ex\h{${\scriptscriptstyle\bullet}$}}
\newcommand{\ssc}{\,\raise.15ex\h{${\scriptstyle\circ}$}\,}
\newcommand{\ds}{\rlap{\raise.4ex\h{$\downarrow$}}{\downarrow}}
\newcommand{\tauZ}{\tau_Z}
\newcommand{\C}{{\mathbb C}}
\newcommand{\D}{{\mathbb D}}
\newcommand{\N}{{\mathbb N}}
\newcommand{\PP}{{\mathbb P}}
\newcommand{\Q}{{\mathbb Q}}
\newcommand{\R}{{\mathbb R}}
\newcommand{\Z}{{\mathbb Z}}
\newcommand{\B}{{\mathcal B}_g}
\newcommand{\BB}{\widetilde{\mathcal B}_g}
\newcommand{\CC}{{\mathcal C}_g}
\newcommand{\CCt}{\widetilde{\mathcal C}_g}
\newcommand{\Dc}{{\mathcal D}}

\newcommand{\G}{{\mathcal G}}
\newcommand{\LL}{{\mathcal L}}
\newcommand{\OO}{{\mathcal O}}
\newcommand{\Rc}{{\mathcal R}}
\newcommand{\Ct}{\widetilde{C}}

\newcommand{\E}{\widetilde{E}}
\newcommand{\HH}{\widetilde{H}}
\newcommand{\I}{\widetilde{I}}

\newcommand{\Hdf}{H_{\df\wedge}}
\newcommand{\X}{\widetilde{X}}
\newcommand{\x}{\widetilde{x}}
\newcommand{\dd}{{\partial}}
\newcommand{\ddd}{{\rm d}}
\newcommand{\dti}{\dd_t^{-1}}
\newcommand{\df}{{\rm d}f}
\newcommand{\ddh}{{\rm d}h}
\newcommand{\sKf}{{}^s\!K^{\ssb}_f}
\newcommand{\spKf}{{}^{s}{}'\!K_{f'}^{\ssb}}
\newcommand{\la}{\lambda}
\newcommand{\mm}{{\mathfrak m}}
\newcommand{\mmm}{\overline{\mm}}
\newcommand{\Om}{\Omega}
\newcommand{\Si}{\Sigma}
\newcommand{\Sg}{{\rm Sing}}
\newcommand{\Sh}{{\mathcal S}\hskip-.5pt h}
\newcommand{\Sp}{{\rm Sp}}
\newcommand{\M}{{}\,\overline{\!M}{}}
\newcommand{\RR}{{}\,\overline{\!R}{}}
\newcommand{\Yt}{\widetilde{Y}}
\newcommand{\Zt}{\widetilde{Z}}
\newcommand{\Gr}{{\rm Gr}}
\newcommand{\Ext}{{\rm Ext}}
\newcommand{\Hom}{{\rm Hom}}
\newcommand{\bl}{\bigl}
\newcommand{\br}{\bigr}
\newcommand{\into}{\hookrightarrow}
\newcommand{\onto}{\mathop{\rlap{$\to$}\hskip2pt\hbox{$\to$}}}
\newcommand{\simto}{\buildrel\sim\over\longrightarrow}
\newcommand{\ges}{\geqslant}
\newcommand{\les}{\leqslant}
\begin{document}
\title[Koszul complexes and spectra]
{Koszul complexes and spectra of projective hypersurfaces with isolated singularities}
\dedicatory{To the memory of Egbert Brieskorn}
\author{Alexandru Dimca}
\address{Univ. Nice Sophia Antipolis, CNRS, LJAD, UMR 7351, 06100 Nice, France.}
\email{dimca@unice.fr}
\author{Morihiko Saito}
\address{RIMS Kyoto University, Kyoto 606-8502 Japan}
\email{msaito@kurims.kyoto-u.ac.jp}
\begin{abstract} For a projective hypersurface $Z$ with isolated singularities, we generalize some well-known assertions in the nonsingular case due to Griffiths, Scherk, Steenbrink, Varchenko, and others about the relations between the Steenbrink spectrum, the Poincar\'e polynomial of the Jacobian ring, and the roots of Bernstein-Sato polynomial for a defining polynomial $f$ up to sign forgetting the multiplicities. We have to use the pole order spectrum and the alternating sum of the Poincar\'e series of certain subquotients of the Koszul cohomologies, and study the pole order spectral sequence. We show sufficient conditions for vanishing or non-vanishing of the differential $d_1$ of the spectral sequence, which are useful in many applications. We prove also symmetries of the dimensions of the subquotients of Koszul cohomologies, which are crucial for computing the roots of BS polynomials. We can deduce that the roots of BS polynomial  whose absolute values are larger than $n-1-n/d$ are determined by the ``torsion part" of the Jacobian ring (modulo the roots of BS polynomial for $Z$) if all the singularities of $Z$ are weighted homogeneous. Here $d=\deg f$ and $n$ is the dimension of the ambient affine space.
\end{abstract}
\maketitle
\centerline{\bf Introduction}
\bsn
Let $f$ be a homogeneous polynomial in the graded $\C$-algebra $R:=\C[x_1,\dots,x_n]$ where $\deg x_i=1$ and $n\ges 2$. Set $d=\deg f$. Consider the shifted Koszul complex
$$\sKf:=K_f^{\ssb}[n]\q\h{with}\q K_f^{\ssb}=(\Om^{\ssb},\df\wedge).$$
\par\nin Here $\Om^j:=\Gamma(\C^n,\Om_{\C^n}^j)$ with $\Om_{\C^n}^j$ algebraic so that the $\Om^j$ are finite free graded $R$-modules, and the degree of $\Om^j$ in $\sKf$ is shifted so that
$${}^s\!K^j_f=\Om^{j+n}(jd)\,\,\,\h{(that is,}\,\,\,{}^s\!K^j_{f,k}=\Om^{j+n}_{jd+k})\,\,\,\h{for}\,\,\,j\in\Z.$$
\par\nin In general the shift of degree by $p$ of a graded module $M$ will denoted by $M(p)$, where the latter is defined by $M(p)_k=M_{k+p}$. Since the dualizing complex for complexes of $R$-modules is given by $\Om^n[n]$, we have the self-duality
$$\D(\sKf):=\R\Hom_R(\sKf,\Om^n[n])=\sKf(nd).$$
\par\nin \sk
In this paper we assume
$$\dim\Sg\,f^{-1}(0)\les 1,\,\,\,\h{and $f$ is not a polynomial of $n{-}1$ variables.}
\leqno({\rm A})$$
\par\nin It is well known, and is easy to show (see for instance Remark~\hl{R1.9}{1.9}\,(iv) below) that this implies
$$H^j(\sKf)=0\q\h{if}\,\,\,j\ne -1,0.$$
\par\nin Define
$$M:=H^0(\sKf),\q N:=H^{-1}(\sKf).$$
\par\nin Let $\mm=(x_1,\dots,x_n)\subset R$, the maximal graded ideal. Set
$$M':=H_{\mm}^0M,\q M'':=M/M'.$$
\par\nin These are finitely generated graded $R$-modules having the decompositions $M=\bigoplus_k M_k$, etc.
In the isolated singularity case we have $M''=N=0$, and $M=M'$.
Generalizing a well-known assertion in the isolated singularity case, one may conjecture that the canonical morphism from $M'$ to the graded quotient of the pole order filtration on the Gauss-Manin system is injective, see Proposition~\hl{P3.5}{3.5} below for a partial evidence.
This is closely related to Question~2 and Remark~\hl{R5.9}{5.9} below.
\sk
Let $y:=\msum_{i=1}^n\,c_ix_i$ with $c_i\in\C$ sufficiently general so that $\{y=0\}\subset\C^n$ is transversal to any irreducible component of $\Sg\,f^{-1}(0)$. Then $M'$ is the $y$-torsion subgroup of $M$, and $M''$, $N$ are finitely generated free graded $\C[y]$-modules of rank $\tauZ$, where $\tauZ$ is the total Tjurina number as in (\hl{4}{4}) below.
Note that there is a shift of the grading on $N$ by $d$ between this paper and \cite{DiSt1}, \cite{DiSt2}.
\sk
Define the (higher) dual graded $R$-modules by
$$D_i(M):=\Ext_R^{n-i}(M,\Om^n)\q(i\in\Z),$$
\par\nin and similarly for $D_i(N)$, etc. From the above self-duality of the Koszul complex $\sKf$, we can deduce the following duality (which is known to the specialists at least by forgetting the grading, see \cite{Pe}, \cite{vStWa}, \cite{Se}, \cite{EyMe}):
\par\htt{T1}{}\msn
{\bf Theorem~1.} {\it There are canonical isomorphisms of graded $R$-modules
\htt{1}{}
$$\aligned D_0(M')=D_0(M)&=M'(nd),\\ D_1(M'')=D_1(M)&=N(nd),\\ D_1(N)&=M''(nd),\endaligned
\leqno(1)$$
\par\nin and $D_i(M)$, $D_i(M')$, $D_i(M'')$, $D_i(N)$ vanish for other $i$.}
\ms
This generalizes a well-known assertion in the isolated singularity case with $M''\,{=}\,N\,{=}\,0$.
Theorem~\hl{T1}{1} implies that the graded $R$-modules $M'$, $M''$, and $N$ are Cohen-Macaulay with dimension $0$, $1$ and $1$ respectively (but $M$ itself is not Cohen-Macaulay). Moreover $M'$ is graded self-dual, and $M''$ and $N$ are graded dual of each other, up to a shift of grading.
\sk
For $k\in\Z$, set
$$\mu_k=\dim M_k,\q\mu'_k=\dim M'_k,\q\mu''_k=\dim M''_k,\q
\nu_k=\dim N_k.$$
\par\nin Let $g:=\msum_{i=1}^nx_i^d$, and $\gamma_k:=\dim(H^nK_g^{\ssb})_k=\dim\bl(\Om^n/\,\msum_{i=1}^n\, x_i^{d-1}\Om^n\br){}_k$, so that
\htt{2}{}
$$\msum_k\,\gamma_k\,t^k=(t^d-t)^n/(t-1)^n.
\leqno(2)$$
\par\nin (Here $g$ can be any homogeneous polynomial of degree $d$ with an isolated singular point.)
It is known (see \cite{Di2}, \cite{DiSt1}, \cite{DiSt2}) that
\htt{3}{}
$$\mu_k=\mu'_k+\mu''_k=\nu_k+\gamma_k\q(k\in\Z),
\leqno(3)$$
\par\nin since the Euler characteristic of a bounded complex is independent of its differential if the components of the complex are finite dimensional.
\sk
By the first assertion of (\hl{1}{1}) together with (\hl{1.1.4}{1.1.4}) for $i=1$ and by (\hl{2}{2}), we get the following symmetries:
\par\htt{C1}{}\msn
{\bf Corollary~1.}\q $\mu'_k=\mu'_{nd-k}\q(k\in\Z)$.
\ms
This is compatible with the symmetry $\gamma_k=\gamma_{nd-k}$. Set $Z:=\{f=0\}\subset Y:=\PP^{n-1}$, and $\Si:=\Sg\,Z$. The total Tjurina number $\tauZ$ is defined by
\htt{4}{}
$$\tauZ:=\msum_{z\in\Si}\,\tau_z\q\h{with}\q\tau_z:=\dim_{\C}\OO_{Y,z}/(h_z,\dd h_z),
\leqno(4)$$
\par\nin where $h_z$ is a local defining equation of $Z$ at $z$, and $\dd h_z$ is the Jacobian ideal of $h_z$ generated by its partial derivatives. By Theorem~\hl{T1}{1}, $M''$ and $N$ are Cohen-Macaulay, and are dual of each other up to a shift of grading. Combining this with the graded local duality (\hl{1.1.4}{1.1.4}) for $i=1$ (see \cite{BrHe}, \cite{Ei}, etc.) together with (\hl{1.9.3}{1.9.3}) below, we get the following (which does not seem to be stated explicitly in the literature). 
\par\htt{C2}{}\msn
{\bf Corollary~2.}\q $\mu''_k+\nu_{nd-k}=\tauZ,\q\delta''_k=\delta''_{(n-1)d-k}\,\,\,$ with $\,\,\,\delta''_k:=\mu''_k-\nu_{k+d}\q(k\in\Z)$.
\ms
The $\delta''_k$ are very important for computations of the roots of Bernstein-Sato polynomials, see \cite{wh}. The calculation of the local cohomology in the local duality is not so trivial (see Remark~\hl{R1.7}{1.7} below), and we can use also an exact sequence as in \cite[Prop.~2.1]{Sch}, see also \cite[Prop.~2.1.5]{Gro}, \cite{SaSch}. Corollary~\hl{C2}{2} can be deduced also from \cite[Theorem 3.1]{Di2}, see Remark~\hl{R1.9}{1.9}\,(i) below. By Corollaries~1 and 2 together with (\hl{3}{3}), we get the following.
\par\htt{C3}{}\msn
{\bf Corollary~3.}\q $\mu'_k=\mu_k+\mu_{nd-k}-\gamma_k-\tauZ,\q
\mu''_k=\tauZ-\mu_{nd-k}+\gamma_k\q(k\in\Z)$.
\ms
This means that $\mu'_k$ and $\mu''_k$ are essentially determined by $\mu_k$ and $\mu_{nd-k}$. Note that $\{\mu''_k\}$ and $\{\nu_k\}$ are weakly increasing sequences of non-negative integers. It is shown that $\{\mu'_k\}$ is log-concave in a certain case, see \cite{Sti}. Assuming $\Sg\,Z\ne\emptyset$, we have $\mu''_k=\nu_k=\tauZ>0$ for $k\gg 0$, hence $M'',N$ are nonzero, although $M'$ may vanish, see Remark~\hl{R1.9}{1.9}\,(iii) below. By Corollary~\hl{C2}{2} and (\hl{3}{3}) we get the following.
\par\htt{C4}{}\msn
{\bf Corollary~4.}\q$\gamma_k-\mu'_k=\mu''_k+\mu''_{nd-k}-\tauZ\q(k\in\Z)$.
\ms
Here a fundamental question seems to be the following.
\msn
{\bf Question~1.} Are both sides of the above equality non-negative?
\msn
This seems to be closely related to the subject treated in \cite{ChDi}, \cite{Di2}, \cite{DiSt1}, \cite{DiSt2}, etc.
We have a positive answer to Question~1 if $n=3$ and $\Si$ is a complete intersection in $\PP^2$ (see \cite{Sti}) or if $f$ has type (I), where $f$ is called type (I) if the following condition is satisfied (and type (II) otherwise):
\htt{5}{}
$$\mu''_k=\tauZ\,\,\,\h{for}\,\,\,k\ges nd/2,\q\h{that is,}\q
\nu_k=0\,\,\,\h{for}\,\,\,k\les nd/2.
\leqno(5)$$
\par\nin By the definition of $N$, the last condition in (\hl{5}{5}) cannot hold if there is a nontrivial relation of very low degree between the partial derivatives of $f$; for instance, in case $f$ is a polynomial of $n-1$ variables (or close to it), see Remark~\hl{R2.9}{2.9} below. However, it holds in relatively simple cases, including the nodal hypersurface case by \cite[Thm.~2.1]{DiSt1}, see Remark~\hl{R2.10}{2.10} below.
\sk
In the type (I) case, we get the $\mu'_k$ by restricting to $k\les nd/2$ (where $\mu'_k+\mu''_k=\mu_k=\gamma_k$ holds) if we know the $\mu''_k$. This can be done for instance in the following case.
\par\htt{P1}{}\msn
{\bf Proposition~1.} {\it Assume $Z$ has only ordinary double points $z_1,\dots,z_{\tauZ}$, and moreover the $z_i$ correspond to linearly independent vectors in $\C^n$ so that $\tauZ=r\les n$. Then
$$\aligned\mu''_k&=\begin{cases}\,0&\h{if}\,\,\,\,\,k<n,\\
\,1&\h{if}\,\,\,\,\,k=n,\\
\tauZ&\h{if}\,\,\,\,\,k>n,\end{cases}\q\q\q
\nu_k=\begin{cases}\,0&\h{if}\,\,\,\,\,k<n(d{-}1),\\
\tauZ-1&\h{if}\,\,\,\,\,k=n(d{-}1),\\
\tauZ&\h{if}\,\,\,\,\,k>n(d{-}1),\end{cases}\\
\mu'_k&=\begin{cases}\,0&\h{if}\,\,\,\,\,k\notin\bl(n,n(d{-}1)\br),\\
\gamma_k-\tauZ&\h{if}\,\,\,\,\,k\in\bl(n,n(d{-}1)\br),\end{cases}\endaligned$$
\par\nin where $\bl(n,n(d{-}1)\br)\subset\R$ denotes an open interval.}
\ms
This follows from Lemma~\hl{L2.1}{2.1} below together with Corollary~\hl{C2}{2} and (\hl{3}{3}). It can also be deduced from the results in \cite{Di2}, and seems to be closely related to \cite[Thm.~2]{DiSaWo}. The situation becomes, however, rather complicated if the number of singular points is large, see \cite{ChDi}, \cite{Di2}, \cite{DiSt1}, \cite{DiSt2}.
\sk
Let $\Sp(f)=\sum_{\alpha}n_{f,\alpha}\,t^{\alpha}\in\Q[t^{1/d}]$ be the {\it Steenbrink spectrum} of $f$ (see \cite{St2}, \cite{St3}) which is normalized as in \cite{St2}. To study the relation with the Koszul cohomologies $M$, $N$ by generalizing the well-known assertion in the isolated singularity case where $M''=N=0$ and $M=M'$ (see \cite{St1} and also \cite{Gri}, \cite{ScSt}, \cite{Va}, etc.), we have to introduce the {\it pole order spectrum} $\Sp_P(f)$ by replacing the Hodge filtration $F$ with the pole order filtration $P$ in \cite{Di1}, \cite{Di3}, \cite{DiSa2}, \cite{DiSt2}. There are certain shifts of the exponents coming from the difference between $F$ and $P$. Here we have the inclusion $F\subset P$ in general, and the equality holds in certain cases (see \cite{Di3}).
We can calculate these spectra explicitly in the case $n=2$, see Propositions~\hl{P3.3}{3.3} and \hl{P3.4}{3.4}.
The relation between the two spectra is, however, quite nontrivial in general (see for instance Example~\hl{E3.7}{3.7} below).
\sk
The reason for which we introduce $\Sp_P(f)$ is that it is related to the Poincar\'e series of $M$, $N$ as follows:
The differential of the de Rham complex $(\Om^{\ssb},\ddd)$ induces a morphism of graded $\C$-vector spaces of degree $-d:$
$$\ddd^{(1)}:N\to M,$$
\par\nin that is, preserving the degree up to the shift by $-d$. Let $H^nA_f^{\ssb}$ denote the Brieskorn module \cite{Bri} (in a generalized sense) which is a graded $\C$-module endowed with actions of $t$, $\dti$, and $t\dd_t$, see \hl{4.2}{4.2} below. Let $(H^nA_f^{\ssb})_{\rm tor}$ be its $t$-torsion (or equivalently, $\dti$-torsion) subspace. It has the kernel filtration $K_{\ssb}$ defined by
\htt{6}{}
$$K_i(H^nA_f^{\ssb})_{\rm tor}:={\rm Ker}^n\,t^i\subset(H^nA_f^{\ssb})_{\rm tor}\q(i\ges 0),
\leqno(6)$$
\par\nin where ${\rm Ker}^n\,t^i$ means that the kernel is taken in $(H^nA_f^{\ssb})_{\rm tor}$.
\sk
One of the main theorems of this paper is as follows:
\par\htt{T2}{}\msn
{\bf Theorem~2.} {\it There are inductively defined morphisms of graded $\C$-vector spaces of degree $-rd:$
$$\ddd^{(r)}:N^{(r)}\to M^{(r)}\q(r\ges 2),$$
\par\nin such that $N^{(r)}$, $M^{(r)}$ are the kernel and the cokernel of $\ddd^{(r-1)}$ respectively, and are independent of $r\gg 0$ $($that is, $\ddd^{(r)}=0$ for $r\gg 0)$, and we have
\htt{7}{}
$$\Sp_P(f)=S(M^{(r)})(t^{1/d})-S(N^{(r)})(t^{1/d})\q(r\gg 0),
\leqno(7)$$
\par\nin where $S(M^{(r)})(t)$, $S(N^{(r)})(t)$ denote the Poincar\'e series of $M^{(r)}$, $N^{(r)}$ for $r\ges 2$.
\sk
Moreover, there are canonical isomorphisms
\htt{8}{}
$${\rm Im}\,\ddd^{(r)}=\Gr^K_{r-1}({\rm Coker}^n\,t)\q(r\ges 2),
\leqno(8)$$
\par\nin where $K_{\ssb}$ is the kernel filtration on $(H^nA_f^{\ssb})_{\rm tor}$ in $(6)$, and ${\rm Coker}^n\,t$ means the cokernel of the action of $t$ on $(H^nA_f^{\ssb})_{\rm tor}$. In particular, $\ddd^{(r)}$ vanishes for any $r\ges 2$ {\rm (}that is, $M^{(r)}=M^{(2)}$, $N^{(r)}=N^{(2)}$ for any $r\ges 2)$ if and only if $H^nA_f^{\ssb}$ is torsion-free.}
\ms
Note that ${\rm Ker}^n\,t^i$ in (\hl{6}{6}) and ${\rm Coker}^n\,t$ in (\hl{8}{8}) can be replaced respectively with ${\rm Ker}^n\,\dd_t^{-i}$ and ${\rm Coker}^n\,\dti$ by using (\hl{4.2.2}{4.2.2}) below.
For the proof of Theorem~\hl{T2}{2} we use the spectral sequence associated with the pole order filtration on the algebraic microlocal Gauss-Manin complex (see (\hl{4.4.4}{4.4.4}) below), and the morphisms $\ddd^{(r)}$ are induced by the differentials $\ddd_r$ of the spectral sequence. (We can also use the usual Gauss-Manin complex instead of the microlocal one.) The last equivalent two conditions in Theorem~\hl{T2}{2} are further equivalent to the $E_2$-degeneration of the (microlocal) pole order spectral sequence, see Corollary~\hl{C4.7}{4.7} below (and also \cite{vSt}). Moreover $(H^nA_f^{\ssb})_{\rm tor}$ is finite dimensional if and only if $Z$ is analytic-locally defined by a weighted homogeneous polynomial at any singular point, see Theorems~\hl{T5.2}{5.2} and \hl{T5.3}{5.3} below.
(Indeed, the if part in the analytic local setting was shown in the second author's master thesis, see for instance \cite[Thm.~3.2]{BaSa} and also \cite{vSt}.)
Note that Theorem~\hl{T5.3}{5.3} gives rather precise information about the kernel of $\ddd^{(1)}$.
This is a refinement of \cite[Thm.~2.4(ii)]{DiSt2}, and is used in an essential way in \cite{DiSa3}.
Theorem~\hl{T5.3}{5.3} implies a sharp estimate for $\max\{k\,|\,\nu_k=0\}$ when $n=3$, see Corollary~\hl{C5.5}{5.5} below.
This assertion is used in an essential way in \cite{DiSe}, and is generalized to the case $n>3$ in \cite[Theorem~9]{DiSa3} (see \cite{Di4} for another approach to the case $n>3$). For applications of Theorem~\hl{T5.3}{5.3} to determinations of the roots of Bernstein-Sato polynomials, see \cite{bha}, \cite{wh}, \cite{nwh}, \cite{deg}.
\sk
In case $(H^nA_f^{\ssb})_{\rm tor}=0$, we can determine the pole order spectrum if we can calculate the morphism $\ddd^{(1)}:N\to M$, although the latter is not so easy in general unless the last conditions in Theorem~\hl{T5.3}{5.3} are satisfied (see also Remark~\hl{R5.9}{5.9} below). Note that the pole order spectral sequence was studied in \cite{vSt} from a slightly different view point in the (non-graded) analytic local case.
\sk
For the moment there are no examples such that the singularities of $Z$ are weighted homogeneous and $(H^nA_f^{\ssb})_{\rm tor}\ne0$. We have the following.
\msn
{\bf Question~2.} Assume all the singularities of $Z$ are weighted homogeneous. Then, is $H^nA_f^{\ssb}$ torsion-free so that the pole order spectral sequence degenerates at $E_2$ and the equality (\hl{7}{7}) holds with $r=2$?
\ms
We have a positive answer in certain cases; for instance, if $n=2$ or $1$ is not an eigenvalue of $(T_z)^d$ for any $z\in\Sg\,Z$ with $T_z$ the monodromy of a local defining polynomial $h_z$ of $(Z,z)$, see Corollary~\hl{C5.4}{5.4} below for a more general condition.
(Question~2 is recently solved positively in \cite{wh}.)
In the above second case, Theorem~\hl{T5.3}{5.3} actually implies the injectivity of $\ddd^{(1)}:N\to M$ (which is a morphism of degree $-d$), and we get the following.
\par\htt{P2}{}\msn
{\bf Proposition~2.} {\it If $(Z,z)$ is weighted homogeneous and $1$ is not an eigenvalue of $(T_z)^d$ for any $z\in\Sg\,Z$, then $H^nA_f^{\ssb}$ is torsion-free and we have}
$$\Sp_P(f)=S(M)(t^{1/d})-S(N)(t^{1/d}).$$
\par\nin \ms
Here the second condition is satisfied if $1$ is not an eigenvalue of $T_z$ and moreover the order of $T_z$ is prime to $d$ for any $z\in\Sg\,Z$.
The second assumption can be replaced with $H^{n-2}(f^{-1}(1),\C)=0$ by \cite{wh} (which solves Question~2 positively), see Remark~\hl{R5.9}{5.9} below for a picture in the optimal case.
Note that Theorem~\hl{T5.2}{5.2} below implies that the pole order spectral sequence {\it cannot\,} degenerate at $E_2$ if $Z$ has an isolated singularity which is {\it not\,} weighted homogenous.
\sk
Note finally that Theorem~\hl{T1}{1} is useful for the determination of the roots of the Bernstein-Sato polynomial supported at 0, since we get the {\it symmetry\,} of the $\delta''_{k}:=\mu''_k-\nu_{k+d}$ with center $\tfrac{(n-1)d}{2}$ by Corollary~\hl{C2}{2}. Let $\Rc_f$, $\Rc_Z$ be the roots of the Bernstein-Sato polynomials $b_f(s)$ and $b_Z(s)$ of $f$ and $Z$ respectively {\it up to sign.} Using Theorem~\hl{T5.3}{5.3} together with \cite[Theorem 2]{cm-b}, we can deduce the following.
\par\htt{T3}{}\msn
{\bf Theorem~3.} {\it Assume all the singularities of $Z$ are isolated and weighted homogeneous. Let $k$ be an integer strictly larger than $(n{-}1)d{-}n$. Assume $\tfrac{k}{d}\notin\Rc_Z$ if $n\ges4$. Then we have}
\htt{9}{}
$$\tfrac{k}{d}\in\Rc_f\iff M'_k\ne 0.
\leqno(9)$$
\par\nin \sk
Here we do {\it not\,} have to use \cite[Theorems 2 and 3]{wh} showing the $E_2$-degeneration of the pole order spectral sequence and the contribution of nonzero $M'_k$ to $\Rc_f$, see \hl{5.10}{5.10} below. 
It is conjectured that the condition $\tfrac{k}{d}\notin\Rc_Z$ follows from the inequality $k>(n{-}1)d{-}n$. This is valid for $n=3$, see for instance \cite{dFEM}. It is easy to show it if $k\ges(n{-}1)d$ or if all the singularities of $Z$ are homogeneous.
In the case $k\ges 2d$ and $n=3$ or $k=2d{-}2$ and $f^{-1}(0)$ is an essential indecomposable reduced central hyperplane arrangement in $\C^3$, the equivalence (\hl{9}{9}) is shown in \cite{Ba} using a completely different method. Note that $\max\Rc_f<2{-}\tfrac{1}{d}$ in the hyperplane arrangement case, see \cite[Theorem 1]{bha}. Another proof of this theorem for $n=3$ is given in \cite[Corollary 7.3]{DiSt3} using an estimate of Castelnuovo-Mumford regularity \cite[Corollary 3.5]{DIM}, see also \cite{Ba}. In general $\max\Rc_f<n$, $\max\Rc_Z<n{-}1$, see \cite{mic}.
\sk
The first named author was partially supported by Institut Universitaire de France.
The second named author was partially supported by Kakenhi 24540039.
\sk
In Section~\hl{S1}{1} we prove Theorem~\hl{T1}{1} after reviewing graded local duality for the convenience of the reader. In Section~\hl{S2}{2} we explain some methods to calculate the Koszul cohomologies in certain cases. In Section~\hl{S3}{3} we recall some basics from the theory of spectra, and prove Propositions~\hl{P3.3}{3.3}, \hl{P3.4}{3.4}, and \hl{P3.5}{3.5}. In Section~\hl{S4}{4} we prove Theorem~\hl{T2}{2} after reviewing some facts from Gauss-Manin systems and Brieskorn modules. In Section~\hl{S5}{5} we calculate $\ddd^{(1)}$ in certain cases, and prove Theorems~\hl{T5.2}{5.2} and \hl{T5.3}{5.3}.
\bs\bs\htt{S1}{}
\vbox{\centerline{\bf 1. Graded local cohomology and graded duality}
\bsn
In this section we prove Theorem~\hl{T1}{1} after reviewing graded local duality for the convenience of the reader.}
\par\htt{1.1}{}\msn
{\bf 1.1.~Graded local duality.} Let $R=\C[x_1,\dots,x_n]$, and $\mm=(x_1,\dots,x_n)\subset R$. Set
\htt{1.1.1}{}
$$\Om^k=\Gamma(\C^n,\Om_{\C^n}^k)\q(k\in\Z).
\leqno(1.1.1)$$
\par\nin Here $\Om_{\C^n}^k$ is algebraic, and $\Om^k$ is a finite free graded $R$-module with $\deg x_i=\deg dx_i=1$.
\sk
For a bounded complex of finitely generated graded $R$-modules $M^{\ssb}$, define
\htt{1.1.2}{}
$$\aligned&D_i(M^{\ssb}):=\Ext_R^{n-i}(M^{\ssb},\Om^n)= H^{-i}\bl(\D(M^{\ssb})\br)\\ &\h{with}\q\D(M^{\ssb}):=\R\Hom_R(M^{\ssb},\Om^n[n]),\endaligned
\leqno(1.1.2)$$
\par\nin where $\D(M^{\ssb})$ can be defined by taking a graded free resolution $P^{\ssb}\to M^{\ssb}$.
\sk
For a finitely generated graded $R$-module $M$, set
\htt{1.1.3}{}
$$H_{\mm}^0M:=\{a\in M\mid
\mm^ka=0\,\,\,\h{for}\,\,\,k\gg 0\}.
\leqno(1.1.3)$$
\par\nin Let $H^i_{\mm}M$ be the cohomological right derived functors $(i\in\N)$. These are defined by taking a graded injective resolution of $M$. We can calculate them by taking a graded free resolution of $M$ as is explained in textbooks of commutative algebra, see for instance \cite{BrHe}, \cite{Ei}. Indeed, $H^i_{\mm}R=0$ for $i\ne n$, and
$$H^n_{\mm}R=\C[x_1^{-1},\dots,x_n^{-1}]\tfrac{1}{x_1\dots x_n},$$
\par\nin where the right-hand side is identified with a quotient of the graded localization of $R$ by $x_1\cdots x_n$ (using a Cech calculation). We then get the {\it graded local duality} for finitely generated graded $R$-modules $M$:
\htt{1.1.4}{}
$$D_i(M)_k={\rm Hom}_{\C}((H^i_{\mm}M)_{-k},\C)\q(k\in\Z,\,i\ges 0),
\leqno(1.1.4)$$
\par\nin see loc.~cit. (Indeed, this can be reduced to the case $M=R$ by the above argument.)
\msn
{\bf Remarks~1.2.} (i) The functors $H^i_{\mm}$ and $D_i$ are compatible with the corresponding functors for non-graded $R$-modules under the forgetful functor, and moreover, the latter functors are compatible with the corresponding sheaf-theoretic functors as is well known in textbooks of algebraic geometry, see for instance \cite{Ha}. However, the information of the grading is lost by passing to the corresponding sheaf unless we use a sheaf with $\C^*$-action.
\sk
(ii) If $M$ is a finitely generated graded $R$-module, then it is well known that
\htt{1.2.1}{}
$$D_i(M)=0\q\h{for}\,\,\,i<0.
\leqno(1.2.1)$$
\par\nin \msn
{\bf 1.3.~Spectral sequences.} For a bounded complex of finitely generated graded $R$-modules $M^{\ssb}$, we have a spectral sequence
\htt{1.3.1}{}
$$'\!E_2^{p,q}(M^{\ssb})=D_{-p}(H^{-q}M^{\ssb})\Longrightarrow D_{-p-q}(M^{\ssb}).
\leqno(1.3.1)$$
\par\nin This can be defined for instance by taking graded free resolutions of $H^iM^{\ssb}$ and ${\rm Im}\,d^{\,i}$ for $i\in\Z$, and then extending these to a graded free resolution of $M^{\ssb}$ by using the short exact sequences
$$0\to{\rm Im}\,d^{\,i-1}\to{\rm Ker}\,d^{\,i}\to H^iM^{\ssb}\to 0,\q 0\to{\rm Ker}\,d^{\,i}\to M^i\to{\rm Im}\,d^{\,i}\to 0,$$
\par\nin as is explained in classical books about spectral sequences. We can also construct (\hl{1.3.1}{1.3.1}) by using the filtration $\tau_{\les -q}$ on $M^{\ssb}$ as in \cite{De}. (Note that $\D(M^{\ssb}[p])=\D(M^{\ssb})[-p]$.)
\sk
Applying (\hl{1.3.1}{1.3.1}) to $\D(M^{\ssb})$ and using $\D(\D(M^{\ssb}))=M^{\ssb}$, we get
\htt{1.3.2}{}
$$''\!E_2^{p,q}(M^{\ssb})=D_{-p}(D_q(M^{\ssb}))\Longrightarrow H^{p+q}M^{\ssb}.
\leqno(1.3.2)$$
\par\nin \par\htt{L1.4}{}\msn
{\bf Lemma~1.4.} {\it Let $\Sh(M)$ denote the coherent sheaf on $X:=\C^n$ corresponding to a finitely generated graded $R$-module $M$. Then we have the following equivalence.}
\htt{1.4.1}{}
$$\aligned H_{\mm}^0M=M&\iff{\rm supp}\,\Sh(M)\subset\{0\}\\ &\iff\h{$M$ is finite dimensional over $\C$},\\ &\iff D_i(M)=0\,\,\,\h{for any}\,\,\,i\ne 0.\endaligned
\leqno(1.4.1)$$
\par\nin \msn
{\it Proof.} This is almost trivial except possibly for the last equivalence. It can be shown by restricting to a sufficiently general point of the support of $\Sh(M)$ in case the support has positive dimension. Here we use the assertion that the dual $\D(\Sh(M))$ is compatible with the direct image under a closed embedding, and this follows from Grothendieck duality for closed embeddings as is well known, see for instance \cite{Ha}. This finishes the proof of Lemma~\hl{L1.4}{1.4}.
\ms
The following is well known, see \cite{BrHe}, \cite{Ei}, etc. We note here a short proof for the convenience of the reader.
\par\htt{P1.5}{}\msn
{\bf Proposition~1.5.} {\it Let $M$ be a finitely generated $R$-module. Set $m:=\dim{\rm supp}\,\Sh(M)$. Then}
\htt{1.5.1}{}
$$D_i(M)=0\,\,\,\h{for}\,\,\,i>m.
\leqno(1.5.1)$$
\par\nin \msn
{\it Proof.} There is a complete intersection $Z$ of dimension $m$ in $X={\rm Spec}\,R$ such that $M$ is annihilated by the ideal $I_Z$ of $Z$, that is, $M$ is an $R_Z$-module with $R_Z:=R/I_Z$, and $I_Z$ is generated by a regular sequence $(g_i)_{i\in[1,n-m]}$ of $R$ with $g_iM=0$. (Here $M$ is not assumed graded.) Set
$$\omega_Z=\Ext_R^{n-m}(R_Z,\Om^n).$$
\par\nin This is called the canonical (or dualizing) module of $Z$. We then get
\htt{1.5.2}{}
$$D_i(M)=\Ext_{R_Z}^{-i}(M,\omega_Z[m]),
\leqno(1.5.2)$$
\par\nin by Grothendieck duality for the closed embedding $i_Z:Z\into X$, see for instance \cite{Ha}, etc. Indeed, taking an injective resolution $G$ of $\Om^n[n]$, one can show (\hl{1.5.2}{1.5.2}) by using the canonical isomorphism
$$\Hom_{R_Z}(M,\Hom_R(R_Z,G))=\Hom_R(M,G).$$
\par\nin Since the right-hand side of (\hl{1.5.2}{1.5.2}) vanishes for $i>m$, the assertion follows.
\par\htt{C1.6}{}\msn
{\bf Corollary~1.6.} {\it Let $M$ be a finitely generated graded $R$-module with $\dim{\rm supp}\,\Sh(M)=1$. Then we have a short exact sequence
\htt{1.6.1}{}
$$0\to D_0(D_0(M))\to M\to D_1(D_1(M))\to 0,
\leqno(1.6.1)$$
\par\nin together with}
\htt{1.6.2}{}
$$D_0(D_1(M))=0,\q D_1(D_0(M))=0.
\leqno(1.6.2)$$
\par\nin \msn
{\it Proof.} By Proposition~\hl{P1.5}{1.5} we get
$$''\!E_2^{p,q}(M)=0\q\h{if}\q(p,q)\notin[-1,0]\times[0,1].$$
\par\nin So the spectral sequence (\hl{1.3.2}{1.3.2}) degenerates at $E_2$ in this case, and the assertion follows.
\par\htt{R1.7}{}\msn
{\bf Remark~1.7.} Let $M$ be a graded $R$-module of dimension 1, that is, $C:={\rm supp}\,\Sh(M)$ is one-dimensional. Let $I_M\subset R$ be the annihilator of $M$. Set $\RR:=R/I_M$. Let $y\in R$ be a general element of degree 1 whose restriction to any irreducible component of $C$ is nonzero. Set $R':=\C[y]\subset R$. Let $\mmm$, $\mm'$ be the maximal graded ideals of $\RR$, $R'$. Let $H^i_{(R,\mm)}M$ denote $H^i_{\mm}M$, and similarly for $H^i_{(\RR,\mmm)}M$, etc. (to avoid any confusion). There are canonical morphisms
$$(R,\mm)\to(\RR,\mmm)\leftarrow(R',\mm'),$$
\par\nin and they imply canonical morphisms
\htt{1.7.1}{}
$$H^i_{(R,\mm)}M\leftarrow H^i_{(\RR,\mmm)}M\to H^i_{(R',\mm')}M.
\leqno(1.7.1)$$
\par\nin Indeed, any graded injective resolution of $M$ over $\RR$ can be viewed as a quasi-isomorphism over $R$ or $R'$, and we can further take its graded injective resolution over $R$ or $R'$, which induces the above morphisms.
\sk
These morphisms are isomorphisms since they are isomorphisms forgetting the grading as is well known. (Note that the morphisms ${\rm Spec}\,R\leftarrow{\rm Spec}\,\RR\to {\rm Spec}\,R'$ are proper. Here it is also possible to use the graded local duality together with Grothendieck duality.) Using the long exact sequence associated with the local cohomology and the localization, we can show
\htt{1.7.2}{}
$$H^1_{(R',\mm')}M=M[y^{-1}]/M.
\leqno(1.7.2)$$
\par\nin So we get the following canonical isomorphism (as graded $R'$-modules):
\htt{1.7.3}{}
$$H^1_{\mm}M=M[y^{-1}]/M.
\leqno(1.7.3)$$
\par\nin This also follows from an exact sequence in \cite[Prop.~2.1]{Sch} (see also \cite[Prop.~2.1.5]{Gro} and \cite{SaSch}, etc.)
\msn
{\bf 1.8.~Proof of Theorem~\hl{T1}{1}.} As is explained in the introduction, we have the self-duality
$$\D(\sKf)=\sKf(nd),$$
\par\nin which implies the isomorphisms of graded $R$-modules
\htt{1.8.1}{}
$$D_i(\sKf)=H^{-i}(\sKf)(nd).
\leqno(1.8.1)$$
\par\nin Consider the spectral sequence (\hl{1.3.1}{1.3.1}) for $M^{\ssb}=\sKf$. By Proposition~\hl{P1.5}{1.5} applied to $M$, $N$, this degenerates at $E_2$. Combining this with (\hl{1.8.1}{1.8.1}), we thus get
\htt{1.8.2}{}
$$D_1(M)=N(nd),\q D_0(N)=0,
\leqno(1.8.2)$$
\par\nin together with a short exact sequence
\htt{1.8.3}{}
$$0\to D_0(M)\to M(nd)\to D_1(N)\to 0.
\leqno(1.8.3)$$
\par\nin By (\hl{1.6.2}{1.6.2}) in Corollary~\hl{C1.6}{1.6} and Proposition~\hl{P1.5}{1.5} applied to $M$, $N$, the proof of Theorem~\hl{T1}{1} is then reduced to showing that (\hl{1.8.3}{1.8.3}) is naturally identified, up to the shift of grading by $nd$, with
\htt{1.8.4}{}
$$0\to M'\to M\to M''\to 0.
\leqno(1.8.4)$$
\par\nin For this, it is enough to show
\htt{1.8.5}{}
$$H^0_{\mm}D_0(M)=D_0(M),\q H^0_{\mm}D_1(N)=0.
\leqno(1.8.5)$$
\par\nin However, the first equality is equivalent to the vanishing of $D_i(D_0(M))$ for $i\ne 0$ by Lemma~\hl{L1.4}{1.4}, and follows from (\hl{1.6.2}{1.6.2}) in Corollary~\hl{C1.6}{1.6} and Proposition~\hl{P1.5}{1.5} applied to $D_0(M)$. The second equality follows for instance from the local duality (\hl{1.1.4}{1.1.4}) for $i=0$ together with (\hl{1.6.2}{1.6.2}) in Corollary~\hl{C1.6}{1.6} applied to $N$. Thus (\hl{1.8.5}{1.8.5}) is proved. This finishes the proof of Theorem~\hl{T1}{1}.
\par\htt{R1.9}{}\msn
{\bf Remarks~1.9.} (i) Corollary~\hl{C2}{2} can be deduced also from \cite[Thm.~3.1]{Di2}. Indeed, by the argument in Section~2 there, we can deduce
\htt{1.9.1}{}
$${\rm def}_{k-n}\Si_f=\tauZ -\mu_k'',
\leqno(1.9.1)$$
\par\nin where ${\rm def}_{k-n}\Si_f$ is as in \cite{Di2}. Moreover, Thm.~3.1 there gives 
\htt{1.9.2}{}
$${\rm def}_{k-n}\Si_f=\mu_{nd-k}-\gamma_{nd-k}=\nu_{nd-k}.
\leqno(1.9.2)$$
\par\nin So Corollary~\hl{C2}{2} follows.
\sk
(ii) It is well known that
\htt{1.9.3}{}
$$\dim_{\C}M''_k=\dim_{\C}M_k=\tauZ\q\h{if}\,\,\,k\gg 0.
\leqno(1.9.3)$$
\par\nin Indeed, the first equality of (\hl{1.9.3}{1.9.3}) is trivial, and it is enough to show the last equality. Changing the coordinates, we may assume $x_n=y$, where $y$ is as in the introduction. On $\{x_n\ne 0\}\subset\C^n$, we have the the coordinates $x'_1,\dots,x'_n$ defined by $x'_j=x_j/x_n$ for $j\ne n$, and $x'_n=x_n$. Using these, we have $f(x)=x_n^dh(x')$, where $x'=(x'_1,\dots,x'_{n-1})$. This implies that the restriction of $\Sh(M)$ to the generic point of an irreducible component of the support of $M$ corresponding to $z\in Z$ has rank $\tau_z$ in the notation of the introduction. So (\hl{1.9.3}{1.9.3}) follows.
\sk
(iii) Assume $\dim\Sg\,f^{-1}(0)=1$, that is, $\Si=\Sg\,Z\ne\emptyset$. Let $(\dd f)\subset R$ denote the Jacobian ideal of $f$ (generated by the partial derivatives $\dd f/\dd x_i$ of $f$). Then the Jacobian ring $R/(\dd f)$ (which is isomorphic to $M$ as a graded $R$-module up to a shift of grading) is a Cohen-Macaulay ring if and only if $M'=0$. Indeed, these are both equivalent to the condition that $M$ is a Cohen-Macaulay $R$-module (since $\tauZ\ne 0$ and hence $M''\ne 0$). Here Grothendieck duality for closed embeddings is used to show the equivalence with the condition that $R/(\dd f)$ is a Cohen-Macaulay ring. Note that $M'$ might vanish in general, for instance if $f$ is as in Example~\hl{E2.7}{2.7} below or even in case $f=xyz$.
\sk
(iv) Assume $\bigcap_{i=1}^mg_i^{-1}(0)
\subset\C^n$ has codimension $\ges r$, where $g_i\in R\,\,\,(i\in[1,m])$. Then there is a regular sequence $(h_j)_{j\in[1,r]}$ of $R$ with $h_j\in V:=\sum_{i=1}^m\C\,g_i$ by increasing induction on $r$ or $m$. This implies the vanishing of the cohomology of the Koszul complex:
$$H^kK^{\ssb}(R;g_1,\dots,g_m)=0\q(k<r),$$
\par\nin by using the $n$-ple complex structure of the Koszul complex as is well known (see Remark~(v) below). Indeed, we can replace the basis $(g_i)$ of the vector space $V$ so that a different expression of the Koszul complex can be obtained. (However, it is not always possible to choose $h_j$ so that $\sum_{i=1}^mRg_i=\sum_{j=1}^rRh_j$ even if $\bigcap_{i=1}^mg_i^{-1}(0)$ has pure codimension $r$ unless $(g_i)$ is already a regular sequence, that is, $r=m$.)
\sk
(v) For $g_i\in R\,\,(i\in[1,m])$, the Koszul complex $K^{\ssb}(R;g_1,\dots,g_m)$ can be identified with the associated single complex of the $m$-ple complex whose $(j_1,\dots,j_m)$-component is $R$ for $(j_1,\dots,j_m)\in[0,1]^m$, and $0$ otherwise, where its $i$-th differential $\ddd_i$ is defined by the multiplication by $g_i$ on $R$.
\sk
(vi) Theorem~\hl{T1}{1} holds with $\df$ in the definition of the Koszul complex replaced by a 1-form $\omega=\sum_{i=1}^n g_idx_i$ if the $g_i$ are homogeneous polynomials of degree $d{-}1$ such that $\bigcap_ig_i^{-1}(0)\subset\C^n$ is at most 1-dimensional. See \cite{Pe}, \cite{vStWa} for the (non-graded) analytic local case.
\bs\bs\htt{S2}{}
\vbox{\centerline{\bf 2. Calculation of the Koszul cohomologies}
\bsn
In this section we explain some methods to calculate the Koszul cohomologies in certain cases.}
\par\htt{L2.1}{}\msn
{\bf Lemma~2.1.} {\it Let $r$ be the dimension of the vector subspace of $\C^n$ generated by the one-dimensional vector subspaces corresponding to the singular points of $Z$. Then}
$$\mu''_n=1,\q\mu''_{n+1}\ges r.$$
\par\nin \msn
{\it Proof.} Let $\Si'$ be a subset of $\Si\,(=\Sg\,Z)$ corresponding to linearly independent $r$ vectors of $\C^n$. Let $I_{\Si'}$ be the (reduced) graded ideal of $R$ corresponding to $\Si'$. There is a canonical surjection
\htt{2.1.1}{}
$$M\to\M:=\Om^n/I_{\Si'}\,\Om^n.
\leqno(2.1.1)$$
\par\nin The target is a free graded $\C[y]$-module of rank $r$, where $y$ is as in the introduction, and it has free homogeneous generators $w_i\,(i\in[1,r])$ with $\deg w_1=n$ and $\deg w_i=n+1$ for $i>1$. So the surjection (\hl{2.1.1}{2.1.1}) factors through $M''$, and the assertion follows.
\par\htt{P2.2}{}\msn
{\bf Proposition~2.2.} {\it Let $f=f_1+f_2$ with $f_1\in\C[x_1,\dots,x_{n_1}]$, $f_2\in\C[x_{n_1+1},\dots,x_n]$, where $1\,{<}\,n_1\,{<}\,n-1$. Assume the dimensions of the singular loci of $f_1^{-1}(0)\,{\subset}\,\C^{n_1}$ and $f_2^{-1}(0)\,{\subset}\,\C^{n-n_1}$ are respectively $1$ and $0$. Then there are isomorphisms of graded $R$-modules
$$M'=M'_{(1)}\sotim_{\C}M'_{(2)},\q M''=M''_{(1)}\sotim_{\C}M'_{(2)},\q N=N_{(1)}\sotim_{\C}M'_{(2)},$$
\par\nin and, setting $S(\mu):=\msum_k\,\mu_k\,t^k\in\Z[[t]]$, etc., we have the equalities
$$S(\mu')=S(\mu'_{(1)})\,S(\mu'_{(2)}),\q S(\mu'')=S(\mu''_{(1)})\,S(\mu'_{(2)}),\q S(\nu)=S(\nu_{(1)})\,S(\mu'_{(2)}),$$
\par\nin where $M'_{(i)}$, $M''_{(i)}$, $N_{(i)}$, and $\mu'_{(i),k}$, $\mu''_{(i),k}$, $\nu_{(i),k}\,\,(k\in\Z)$ are defined for $f_i\,\,\,(i=1,2)$.}
\msn
{\it Proof.} Using the $n$-ple complex structure of the Koszul complex as in Remark~\hl{R1.9}{1.9}\,(v), we get the canonical isomorphism
$$\sKf={}^s\!K^{\ssb}_{f_1}\otimes_{\C}{}^s\!K^{\ssb}_{f_2},$$
\par\nin where ${}^s\!K^{\ssb}_{f_1}$ is defined by using the subring $\C[x_1,\dots,x_{n_1}]$, and similarly for ${}^s\!K^{\ssb}_{f_2}$. Since $f_2^{-1}(0)$ has an isolated singularity, $K^{\ssb}_{f_2}$ is naturally quasi-isomorphic to $M'_{(2)}$. We get hence
$$M=M_{(1)}\sotim_{\C}M'_{(2)},\q N=N_{(1)}\sotim_{\C}M'_{(2)}.$$
\par\nin Moreover, the freeness of $M''_{(1)}\sotim_{\C}M'_{(2)}$ over $\C[y]$ can be shown using an appropriate filtration of $M'_{(2)}$, where $y$ is as in the introduction. These imply that the following two short exact sequences are identified with each other:
$$\aligned 0\to M'\to&M\to M''\to 0,\\ 0\to M'_{(1)}\sotim_{\C}M'_{(2)}\to M_{(1)}&\sotim_{\C}M'_{(2)}\to M''_{(1)}\sotim_{\C}M'_{(2)}\to 0.\endaligned$$
\par\nin So the assertion follows.
\ms
For the proof of Proposition~\hl{P2.6}{2.6} below, we need the following lemma. Essentially this may be viewed as a special case of \cite[Prop.~13]{ChDi}, see Remark~\hl{R2.5}{2.5} below. We note here a short proof of the lemma using Corollaries~1 and 2 and (\hl{3}{3}) in the introduction for the convenience of the reader.
\par\htt{L2.3}{}\msn
{\bf Lemma~2.3.} {\it Assume $n\,{=}\,2$. Let $r$ be the number of irreducible components of $f^{-1}(0)\,{\subset}\,\C^2$. Then $\tauZ=d-r$, and we have for $k\in\Z$
\htt{2.3.1}{}
$$\aligned\mu'_k&=\max(r-1-|d-k|,0),\\ \mu''_k&=(k-1)_{[0,\tauZ]},\\ \nu_k&=(k-d-r+1)_{[0,\tauZ]},\endaligned
\leqno(2.3.1)$$
\par\nin where $x_{[\alpha,\beta]}$ for $x,\alpha,\beta\in\Z$ with $\alpha<\beta$ is defined by}
$$x_{[\alpha,\beta]}=\begin{cases}\alpha&\h{if}\,\,\,x\les\alpha,\\ x&\h{if}\,\,\,\alpha\les x\les\beta,\\ \beta&\h{if}\,\,\,\beta\les x.\end{cases}$$
\par\nin \msn
{\it Proof.} We have the decomposition
$$f=\mprod_{i=1}^r\,g_i^{m_i},$$
\par\nin with $\deg g_i=1$ and $m_i\ges 1$. For $z\in\PP^1$ corresponding to $g_i^{-1}(0)\subset\C^2$, we have
$$\tau_z=m_i-1,\q\h{and hence}\q\tauZ=d-r.$$
\par\nin Setting
$$f':=\mprod_{i=1}^r\,g_i^{m_i-1},$$
\par\nin we get
$$M''=\Om^2/f'\,\Om^2.$$
\par\nin Indeed, the right-hand side is a quotient graded $R$-module of $M$, and is a free graded $\C[y]$-module of rank $\tauZ$. Since $\deg f'=\tauZ$, this implies
$$\mu''_k=(k-1)_{[0,\tauZ]}.$$
\par\nin Using Corollary~\hl{C2}{2}, we then get
$$\nu_k=d-r-(2d-k-1)_{[0,\tauZ]}=(k-d-r+1)_{[0,\tauZ]}.$$
\par\nin Here note that
$$\nu_k=0\q\h{if}\,\,\,k\les d.$$
\par\nin For $n=2$ and $k\les d$, we have
$$\gamma_k=\max(k-1,0).$$
\par\nin By (\hl{3}{3}) in the introduction we then get for $k\les d$
$$\mu'_k=\gamma_k-\mu''_k=\max(k-1-\tauZ,0).$$
\par\nin The formula for $k\ges d$ follows by using the symmetry in Corollary~\hl{C1}{1}. This finishes the proof of Lemma~\hl{L2.3}{2.3}.
\ms
By an easy calculation we see that Lemma~\hl{L2.3}{2.3} is equivalent to the following.
\par\htt{C2.4}{}\msn
{\bf Corollary~2.4.} {\it With the notation and the assumption of Lemma~$(\hl{L2.3}{2.3})$, we have
\htt{2.4.1}{}
$$\aligned S(\mu')&=\Si^{\,[1,r-1]}\,\Si^{\,[d{-}r{+}1,d{-}1]},\\ S(\mu'')&=\Si^{\,[1,\infty)}\,\Si^{\,[1,d-r]},\\ S(\nu)&=\Si^{\,[d+r-1,\infty)}\,\Si^{\,[1,d-r]},\endaligned
\leqno(2.4.1)$$
\par\nin where $S(\mu')$ is as in Proposition~$(\hl{P2.2}{2.2})$, and $\Si^{\,[a,b]}$, $\Si^{\,[a,\infty)}$for $a,b\in\N$ are defined by}
\htt{2.4.2}{}
$$\Si^{\,[a,\infty)}:=\msum_{k\ges a}\,t^k,\q\Si^{\,[a,b]}:=\msum_{k=a}^b\,t^k\,\,\,\h{if}\,\,\,a\les b,
\,\,\,\h{and}\,\,\,0\,\,\,\h{otherwise}.
\leqno(2.4.2)$$
\par\nin \par\htt{R2.5}{}\msn
{\bf Remark~2.5.} With the notation and the assumption of Corollary~\hl{C2.4}{2.4}, the following is shown in \cite[Example~14~(i)]{Di2} as a corollary of Prop.~13 there
\htt{2.5.1}{}
$$S(\mu)=t^2(1-2t^{d-1}+t^{d+r-2})/(1-t)^2.
\leqno(2.5.1)$$
\par\nin By Corollaries~2 and 3 together with (\hl{3}{3}) in the introduction, this is essentially equivalent to the equalities in (\hl{2.4.1}{2.4.1}). Indeed, it seems rather easy to deduce (\hl{2.5.1}{2.5.1}) from (\hl{2.4.1}{2.4.1}). For the converse some calculation seems to be needed. (The details are left to the reader.)
\ms
In case $n_1=2$, we can calculate $\mu'_{(1),k}$, $\mu''_{(1),k}$, $\nu_{(1),k}$ for $f_1$ by Lemma~\hl{L2.3}{2.3}, and get the following.
\par\htt{P2.6}{}\msn
{\bf Proposition~2.6.} {\it Assume $f=f_1+f_2$ as in Proposition~$(\hl{P2.2}{2.2})$ with $n_1=2$. Let $r$ be the number of the irreducible components of $f_1^{-1}(0)\subset\C^2$. Then, under the assumption of Proposition~$(\hl{P2.2}{2.2})$, we have
$$\aligned S(\mu')&=\Si^{\,[1,r-1]}\,\Si^{\,[d-r+1,d-1]}\,(\Si^{\,[1,d-1]})^{n-2},\\ S(\mu'')&=\Si^{\,[1,\infty)}\,\Si^{\,[1,d-r]}\,(\Si^{\,[1,d-1]})^{n-2},\\ S(\nu)&=\Si^{\,[d+r-1,\infty)}\,\Si^{\,[1,d-r]}\,(\Si^{\,[1,d-1]})^{n-2},\endaligned$$
\par\nin where $\Si^{\,[a,b]}$, $\Si^{\,[a,\infty)}$ are as in $(2.4.2)$.}
\msn
{\it Proof.} The assertion follows from Corollary~\hl{C2.4}{2.4} and Proposition~\hl{P2.2}{2.2}, since $S(\mu')$ in the isolated singularity case is invariant by $\mu$-constant deformation, and is given by (\hl{2}{2}) in the introduction.
\par\htt{E2.7}{}\msn
{\bf Example~2.7.} Let $f$ be as in Theorem~\hl{T1}{1}, and assume further
$$f\in\C[x_1,\dots,x_{n-1}]\subset\C[x_1,\dots,x_n].$$
\par\nin Then $f$ has an isolated singularity at the origin of $\C^{n-1}$. Set
$$\gamma'_j:=\dim_{\C}\bl(\Om^{\prime\,n-1}/\df\wedge\Om^{\prime\,n-2}
\br){}_j\q\h{with}\q\Om^{\prime\,k}:=\Gamma(\C^{n-1},\Om_{\C^{n-1}}^k).$$
\par\nin We have the symmetry
\htt{2.7.1}{}
$$\gamma'_j=\gamma'_{(n-1)d-j}.
\leqno(2.7.1)$$
\par\nin In this case, we have $M'=0$, and
\htt{2.7.2}{}
$$\mu_k=\mu''_k=\msum_{j\les k-1}\,\gamma'_j,\q
\nu_k=\msum_{j\les k-d}\,\gamma'_j=\msum_{j\ges nd-k}\,\gamma'_j,
\leqno(2.7.2)$$
\par\nin where the last equality follows from the symmetry (\hl{2.7.1}{2.7.1}), and Corollary~\hl{C2}{2} is verified directly in this case.
\sk
Equivalently, $\mu''_k=\mu_k$ and $\nu_k$ are given as follows:
\htt{2.7.3}{}
$$\aligned S(\mu)&=\Si^{\,[1,\infty)}\,(\Si^{\,[1,d-1]})^{n-1},\\ S(\nu)&=\Si^{\,[d,\infty)}\,(\Si^{\,[1,d-1]})^{n-1},\endaligned
\leqno(2.7.3)$$
\par\nin where $S(\mu)$, etc.\ are as in Proposition~\hl{P2.2}{2.2}, and the order of the coordinates are changed.
\par\htt{E2.8}{}\msn
{\bf Example~2.8.} Assume $n,d\ges 3$. Let
\htt{2.8.1}{}
$$f=x_1^ax_2^{d-a}+\msum_{i=3}^n\,x_i^d\q\h{with}\,\,\,0<a<d.
\leqno(2.8.1)$$
\par\nin We can apply Proposition~\hl{P2.6}{2.6} to this example. More precisely, the calculation of $\mu'_k$, $\mu''_k$ and $\nu_k$ are reduced to the case $n=2$ by Proposition~\hl{P2.2}{2.2}, where $n_1=2$ and
$$f_1=x_1^ax_2^{d-a},\q f_2=\msum_{i=3}^n\,x_i^d.$$
\par\nin The calculation for $f_1$ follows from Lemma~\hl{L2.3}{2.3} or Corollary~\hl{C2.4}{2.4} where $r=2$. For instance, we get in the notation of Proposition~\hl{P2.2}{2.2}
$$\mu'_{(1),k}=\begin{cases}1&\h{if}\,\,\,k=d,\\ 0&\h{if}\,\,\,k\ne d,\end{cases}$$
\par\nin and hence
$$\mu'_k=\mu'_{(2),k+d}=\gamma''_{k+d}\q(k\in\Z),$$
\par\nin where $\gamma''_k$ is as in (\hl{2}{2}) in the introduction with $n$ replaced by $n-2$. By Proposition~\hl{P2.6}{2.6}, we have
\htt{2.8.2}{}
$$\aligned S(\mu')&=t^d\,(\Si^{\,[1,d-1]})^{n-2},\\ S(\mu'')&=\Si^{\,[1,\infty)}\,\Si^{\,[1,d-2]}\,(\Si^{\,[1,d-1]})^{n-2},\\ S(\nu)&=\Si^{\,[d+1,\infty)}\,\Si^{\,[1,d-2]}\,(\Si^{\,[1,d-1]})^{n-2},\endaligned
\leqno(2.8.2)$$
\par\nin where $S(\mu')$, etc.\ are as in Proposition~\hl{P2.2}{2.2}.
\par\htt{R2.9}{}\msn
{\bf Remark~2.9.} If there is a nontrivial relation of degree $k\les d-2$ among the partial derivatives $f_i:=\dd f/\dd x_i$, that is, if there are homogeneous polynomials $g_i$ of degree $k\les d-2$ with $\msum_i\,g_if_i=0$ and $g_i\ne0$ for some $i$, then we have
\htt{2.9.1}{}
$$\nu_{d+n+k-1}\ne 0,
\leqno(2.9.1)$$
\par\nin and hence
\msn
(2.9.2)\q Condition~(\hl{5}{5}) in the introduction does not hold if $(n-2)(d-2)\ges 2(k+1)$.
\msn
Indeed, (\hl{2.9.1}{2.9.1}) follows from the definition $N:=H^{-1}(\sKf)$ since $\deg f_i=d{-}1$.
\sk
This applies to $f$ in Example~\hl{E2.7}{2.7} with $k=0$ since $f_n=0$, and to $f$ in Example~\hl{E2.8}{2.8} with $k=1$ since
$$(d-a)x_1\,f_1=ax_2\,f_2.$$
\par\nin \par\htt{R2.10}{}\msn
{\bf Remark~2.10.} Conditions~(\hl{5}{5}) in the introduction hold in the nodal hypersurface case by \cite[Thm.~2.1]{DiSt1}. Indeed, it is shown there that
\htt{2.10.1}{}
$$\h{$\nu_k=0\,$ if $\,k\les(n_1+1)d\,$ with $\,n\,$ even or $\,k\les(n_1+1)d-1\,$ with $\,n\,$ odd,}
\leqno(2.10.1)$$
\par\nin where $n_1:=[(n{-}1)/2]$. (There is a difference in the grading on $N$ by $d$ between this paper and loc.~cit., and $n$ in this paper is $n+1$ in loc.~cit.)
\bs\bs\htt{S3}{}
\vbox{\centerline{\bf 3. Spectrum}
\bsn
In this section we recall some basics from the theory of spectra, and prove Propositions~\hl{P3.3}{3.3}, \hl{P3.4}{3.4}, and \hl{P3.5}{3.5}.}
\par\htt{3.1}{}\msn
{\bf 3.1.~Hodge and pole order filtrations.} Let $f$ be a homogeneous polynomial of $n$ variables with degree $d$. It is well known that there is a $\C$-local system $L_k$ ($k\in[1,d]$) of rank 1 on $U:=Y\setminus Z$ such that
\htt{3.1.1}{}
$$H^j(U,L_k)=H^j(f^{-1}(1),\C)_{\la}\q\bl(\la=\exp(-2\pi ik/d),
\,k\in[1,d]\br),
\leqno(3.1.1)$$
\par\nin where $H^j(f^{-1}(1),\C)_{\la}$ is the $\la$-eigenspace of the cohomology for the semisimple part of the monodromy, see for instance\ \cite{Di1}, etc. (Note that monodromy in our paper means the one as a local system, see also \cite[Section 1.3]{BuSa}, etc.) Let $\LL_k$ be the meromorphic extension of $L_k\sotim_{\C}\OO_U$. This is a regular holonomic $\Dc_Y$-module, and
\htt{3.1.2}{}
$$H^j\bl(Y,\Om_Y^{\ssb}(\LL_k)\br)=H^j(f^{-1}(1),\C)_{\la}\q\bl(\la=\exp(-2\pi ik/d),
\,k\in[1,d]\br),
\leqno(3.1.2)$$
\par\nin where $\Om_Y^{\ssb}(\LL_k)$ denotes the de Rham complex of $\LL_k$. We have the Hodge and pole order filtrations $F_{\ssb}$ and $P_{\ssb}$ on $\LL_k$ such that
\htt{3.1.3}{}
$$F_i\subset P_i,
\leqno(3.1.3)$$
\par\nin where the equality holds outside the singular points of $Z$, and
$$P_i\LL_k=\begin{cases}\OO_Y(id+k)&\h{if}\,\,\,i\ges 0,\\
\,0&\h{if}\,\,\,i<0,\end{cases}$$
\par\nin see for instance\ \cite[Section 4.8]{cm-b}. (Note that $F$ comes from the Hodge filtration of a mixed Hodge module.) Set $F^i=F_{-i}$, $P^i=P_{-i}$. They induces the Hodge and pole order filtrations on $\Om_Y^{\ssb}(\LL_k)$ such that the $j$-th components of $F^i\,\Om_Y^{\ssb}(\LL_k)$, $P^i\,\Om_Y^{\ssb}(\LL_k)$ are respectively given by
$$\Om_Y^j\sotim_{\OO_Y}F^{i-j}\LL_k,\q
\Om_Y^j\sotim_{\OO_Y}P^{i-j}\LL_k.$$
\par\nin By the isomorphism (\hl{3.1.2}{3.1.2}) they further induce the Hodge and pole order filtrations on the Milnor cohomology $H^j(f^{-1}(1),\C)$. Here $F$ coincides with the Hodge filtration of the canonical mixed Hodge structure. By using the Bott vanishing theorem, $H^{\ssb}\bl(Y,P^i\,\Om_Y^{\ssb}(\LL_k)\br)$ can be calculated by the complex whose $j$-th component is
$$\Gamma(Y,\Om_Y^j\sotim_{\OO_Y}P^{i-j}\LL_k)=
\begin{cases}\Gamma\bl(Y,\Om_Y^j((j-i)d+k)\br)&\h{if}\,\,\,j\ges i,\\ 0&\h{if}\,\,\,j<i.\end{cases}$$
\par\nin But it does not give a strict filtration, and it is not necessarily easy to calculate it.
\sk
Note that the pole order filtration coincides with the one defined by using the Gauss-Manin system, see (\hl{4.4.7}{4.4.7}) and (\hl{4.5.7}{4.5.7}) below.
\par\htt{3.2}{}\msn
{\bf 3.2.~Spectrum.} For $f$ as in \hl{3.1}{3.1}, the spectrum $\Sp(f)=\msum_{\alpha\in\Q}\,n_{f,\alpha}\,t^{\alpha}$ is defined by
\htt{3.2.1}{}
$$\aligned n_{f,\alpha}:=\msum_j\,(-1)^{j-n+1}\dim\Gr^p_F
\HH^j(f^{-1}(1),\C)_{\la}\\
\h{with}\q p=\lfloor n-\alpha\rfloor,\,\,\la=\exp(-2\pi i\alpha),\endaligned
\leqno(3.2.1)$$
\par\nin (see \cite{St2}, \cite{St3}). Here $\HH^j(f^{-1}(1),\C)$ is the reduced cohomology, and we set
\htt{3.2.2}{}
$$\lfloor\alpha\rfloor:=\max\{\,i\in\Z\mid i\les\alpha\,\},\q
\lceil\alpha\rceil:=\min\{\,i\in\Z\mid i\ges\alpha\,\}\q(\alpha\in\R).
\leqno(3.2.2)$$
\par\nin The {\it pole order spectrum} $\Sp_P(f)=\msum_{\alpha\in\Q}\,{}^P\!n_{f,\alpha}t^{\alpha}$ is defined by replacing $F$ with $P$.
\sk
For $j\in\N$, we define $\Sp^j(f)=\msum_{\alpha\in\Q}\,n^j_{f,\alpha}\,t^{\alpha}$ by
\htt{3.2.3}{}
$$\aligned n^j_{f,\alpha}:=\dim\Gr^p_F\HH^{n-1-j}(f^{-1}(1),\C)_{\la}\\
\h{with}\q p=\lfloor n-\alpha\rfloor,\,\,\la=\exp(-2\pi i\alpha),\endaligned
\leqno(3.2.3)$$
\par\nin so that
$$\Sp(f)=\msum_j\,(-1)^j\,\Sp^j(f).$$
\par\nin Similarly $\Sp^j_P(f)=\msum_{\alpha\in\Q}\,{}^P\!n^j_{f,\alpha}\,t^{\alpha}$ is defined by replacing $F$ with $P$.
\sk
Set $Z:=\{f=0\}\subset Y:=\PP^{n-1}$. Let $\pi:(\Yt,\Zt)\to(Y,Z)$ be an embedding resolution, and $E_i$ be the irreducible components of $\Zt$ with $m_i$ the multiplicity of $\Zt$ at the generic point of $E_i$. Let $\alpha=k/d+q\in(0,n)$ with $k\in[1,d]$, $q\in[0,n-1]$. We have by \cite[1.4.3]{BuSa}
\htt{3.2.4}{}
$$n^j_{f,\alpha}=\dim H^{q-j}\bl(\Yt,\Om_{\Yt}^{n-1-q}(\log\Zt)\sotim _{\OO}\,\OO_{\Yt}(-\ell\,\HH+\msum_i\,\lfloor\ell\,m_i/d\rfloor)E_i\br),
\leqno(3.2.4)$$
\par\nin where $\ell:=d-k$, and $\HH$ is the pull-back of a sufficiently general hyperplane $H$ of $Y$.
\sk
In a special case we get the following.
\par\htt{P3.3}{}\msn
{\bf Proposition~3.3.} {\it Assume $n=2$. Set $e:={\rm GCD}(m_i)$ with $m_i$ the multiplicities of the irreducible factors of $f$. Then, for $\alpha=k/d+q\in(0,2)$ with $k\in[1,d]$, $q=0,1$, we have
\htt{3.3.1}{}
$$n^j_{f,\alpha}=\begin{cases}r-1+k-\msum_i\,\lceil km_i/d\rceil&\h{if}\,\,\,j=0,\,q=0,\\
\max\bl(-k-1+\msum_i\,\lceil km_i/d\rceil,\,0\,\br)&\h{if}\,\,\,j=0,\,q=1,\\
1&\h{if}\,\,\,j=1,\,q=1,\,ke/d\in\Z,\\
0&\h{otherwise},\end{cases}
\leqno(3.3.1)$$
\par\nin where $\lceil\alpha\rceil$ is as in $(3.2.2)$.}
\msn
{\it Proof.} We have $\Om_{\PP^1}^1(\log Z)=\OO_{\PP^1}(r-2)$ with $\Yt=Y=\PP^1$, $\Zt=Z$, $\HH=H$, Hence (\hl{3.2.4}{3.2.4}) in this case becomes
$$n^j_{f,\alpha}=\begin{cases}\dim H^0\bl(\PP^1,\Om_{\PP^1}^1(\log Z) (-\ell+\msum_i\,\lfloor\ell\,m_i/d\rfloor)\br)&\h{if}\,\,\,j=0,\,q=0,\\
\dim H^1\bl(\PP^1,\OO_{\PP^1}(-\ell+\msum_i\,\lfloor\ell\,m_i/d\rfloor)\br)&\h{if}\,\,\,j=0,\,q=1,\\
\dim H^0\bl(\PP^1,\OO_{\PP^1}(-\ell+\msum_i\,\lfloor\ell\,m_i/d\rfloor)\br)&\h{if}\,\,\,j=1,\,q=1,\\
0&\h{otherwise},\end{cases}$$
\par\nin and then
\htt{3.3.2}{}
$$n^j_{f,\alpha}=\begin{cases}r-1-\ell+\msum_i\,\lfloor\ell\,m_i/d\rfloor&\h{if}\,\,\,j=0,\,q=0,\\
\max\bl(\ell-1-\msum_i\,\lfloor\ell\,m_i/d\rfloor,\,0\,\br)&\h{if}\,\,\,j=0,\,q=1,\\
\max\bl(-\ell+1+\msum_i\,\lfloor\ell\,m_i/d\rfloor,\,0\,\br)&\h{if}\,\,\,j=1,\,q=1,\\
0&\h{otherwise}.\end{cases}
\leqno(3.3.2)$$
\par\nin Since $\msum_i\,m_i=d$ and $e={\rm GCD}(m_i)$, we have
$$\ell>\msum_i\,\lfloor\ell m_i/d\rfloor\iff\ell m_i/d\notin\Z\,\,(\exists\,i)\iff\ell e/d\notin\Z.$$
\par\nin So (\hl{3.3.1}{3.3.1}) follows (since $\ell=d-k$). This finishes the proof of Proposition~\hl{P3.3}{3.3}.
\ms
We note here an application of Theorem~\hl{T2}{2}, Theorem~\hl{T5.3}{5.3} and Corollary~\hl{C5.4}{5.4})below.
(This will not be used in their proofs.)
\par\htt{P3.4}{}\msn
{\bf Proposition~3.4.} {\it Assume $n=2$. Then $\Sp_P(f)=\Sp_P^0(f)-\Sp_P^1(f)$ is given by
\htt{3.4.1}{}
$$\Sp^j_P(f)=\begin{cases}\msum_k\,(\mu_k-\nu_{k+d})\,t^{k/d}+\bl(\,t^{1/e}+\cdots+t^{(e-1)/e}\,\br)&\h{if}\,\,\,j=0,\\ t\,(\,t^{1/e}+\cdots+t^{(e-1)/e}\,)&\h{if}\,\,\,j=1,\end{cases}
\leqno(3.4.1)$$
\par\nin with $\mu_k$, $\nu_k$ explicitly expressed in Lemma~$(\hl{L2.3}{2.3})$, and $e={\rm GCD}(m_i)$ as in Proposition~{\rm\hl{P3.3}{3.3}}.}
\msn
{\it Proof.} The pole order spectral sequence degenerates at $E_2$ by Corollary~\hl{C5.4}{5.4} below. So the assertion is shown in the case $e=1$, since the last condition implies that $\nu_k^{(2)}=0$. In the general case it is well known that
\htt{3.4.2}{}
$$\HH^0(f^{-1}(1),\C)_{\la}=\begin{cases}\C&\h{if}\,\,\,\la^e=1\,\,\,\h{with}\,\,\,\la\ne 1,\\ \,0&\h{otherwise}\end{cases}.
\leqno(3.4.2)$$
\par\nin By using Theorem~\hl{T5.3}{5.3} and Lemma~\hl{L2.3}{2.3}, this implies
\htt{3.4.3}{}
$$N_{k+d}^{(2)}=\begin{cases}\C&\h{if}\,\,\,k=i\,(d/e)\,\,\h{with}\,\,\,i\in\{1,\dots,e-1\},\\ \,0&\h{otherwise},\end{cases}\leqno(3.4.3)$$
\par\nin where $N^{(2)}\subset N$ is the kernel of $\ddd^{(1)}$. This gives also the information of the coimage of $\ddd^{(1)}$ which is a morphism of degree $-d$. So the correction terms for $\Sp_P^0(f)$ and $\Sp_P^1(f)$ coming form the non-vanishing of $\ddd^{(1)}$ are given respectively by
$$t^{1/e}+\cdots+t^{(e-1)/e}\q\h{and}\q t\,(\,t^{1/e}+\cdots+t^{(e-1)/e}\,).$$
\par\nin So (\hl{3.4.1}{3.4.1}) follows. This finishes the proof of Proposition~\hl{P3.4}{3.4}.
\par\htt{P3.5}{}\msn
{\bf Proposition~3.5.} {\it Assume $f=f_1+f_2$ as in Proposition~$(\hl{P2.2}{2.2})$ with $n_1=2$. Then, under the assumption of Proposition~$(\hl{P2.2}{2.2})$, we have
\htt{3.5.1}{}
$$\mu'_k\les n^0_{f,k/d},\q\mu'_k\les {}^P\!n^0_{f,k/d}\q(k\in\Z),
\leqno(3.5.1)$$
\par\nin where $n^0_{f,k/d}$, ${}^P\!n^0_{f,k/d}$ are as in $(3.2)$.}
\msn
{\it Proof.} The Thom-Sebastiani type theorem holds for $\Sp^0(f)$, $\Sp^0_P(f)$ under the assumption of Proposition~\hl{P2.2}{2.2}, see \hl{4.9}{4.9} below. So the assertion is reduced to the case $f=f_1$ with $n=2$.
The assertion for $\Sp^0_P(f)$ then follows from Proposition~\hl{P3.4}{3.4} and Lemma~\hl{L2.3}{2.3}, where we may assume $r\ges 2$ since $\mu'_k=0$ otherwise.
By using Lemma~\hl{L2.3}{2.3} and Proposition~\hl{P3.3}{3.3} (more precisely, (\hl{3.3.2}{3.3.2}) for $q=0$ and (\hl{3.3.1}{3.3.1}) for $q=1$), the assertion for $\Sp^0(f)$ is reduced to the following trivial inequalities
$$\aligned r-1-(d-k)\les r-1-\ell+\msum_{i=1}^r\,\lfloor\ell\,m_i/d\rfloor&\q(\ell\in[0,d{-}1],\,q=0),\\
r-1+d-(k+d)\les -k-1+\msum_{i=1}^r\,\lceil km_i/d\rceil&\q(k\in[1,d{-}1],\,q=1),\endaligned$$
\par\nin where $\ell=d-k$. (Note that $k$ in Lemma~\hl{L2.3}{2.3} is $k+d$ in the case $q=1$.) This finishes the proof of Proposition~\hl{P3.5}{3.5}.
\par\htt{R3.6}{}\msn
{\bf Remarks~3.6.} (i) If $f$ has an isolated singularity, the equality holds in (\hl{3.5.1}{3.5.1}), and $S(\mu')$ (with $t$ replaced by $t^{1/d}$) coincides with the spectrum $\Sp(f)$, see \cite{St1} and also \cite{Gri}, \cite{ScSt}, \cite{Va}, etc.
It would be interesting if (\hl{3.5.1}{3.5.1}) holds in a more general case.
\ms
(ii) Let $f$ be as in \hl{3.1}{3.1}. Assume $Z\subset\PP^{n-1}$ has only isolated singularities. Let $\alpha'_f$ be the minimal of the exponents of the spectrum for all the singularities of $Z$ (see also Corollary~\hl{C5.5}{5.5} below). Then the multiplicity $n_{f,\alpha}$ of the spectrum $\Sp(f)$ for $\alpha=p/d<\min(\alpha'_f,1)$ can be given by
$$n_{f,p/d}=\binom{p-1}{n-1}\q(p/d<\min(\alpha'_f,1)).$$
\par\nin This follows from a formula for multiplier ideals \cite[Prop.~1]{cm-b} together with \cite{Bu} (see also a remark before \cite[Cor.~1]{cm-b}). This equality holds also for the pole order spectrum since $\mu_p$ is at most the right-hand side of the equality and $F\subset P$ (and $\nu_p=0$ for $p<d$).
\par\htt{E3.7}{}\msn
{\bf Example~3.7.} Let $f=(x^m+y^m)\,x^my^m$ ($m\ges 2$), where $d=3m$, $r=m+2$, $\tauZ=2m-2$.
For $\alpha=k/3m+q$ with $k\in[1,3m]$, $q=0,1$, we have by Proposition~\hl{P3.3}{3.3}
\htt{3.7.1}{}
$$n_{f,\alpha}=\begin{cases}k+1-2\lceil k/3\rceil&\h{if}\,\,\,\alpha\in(0,1],\,\,\,q=0,\\
m-k-1+2\lceil k/3\rceil&\h{if}\,\,\,\alpha\in(1,2),\,\,\,q=1.\end{cases}
\leqno(3.7.1)$$
\par\nin Indeed, $m_i=1$ for $i\in[1,m]$, and $m_i=m$ for $i=m+1,m+2$. Here $e=1$ in the notation of Proposition~\hl{P3.3}{3.3}, and hence $\Sp^1(f)=0$, $\Sp(f)=\Sp^0(f)$ (similarly for $\Sp_P(f)$).
\sk
On the other hand, Lemma~\hl{L2.3}{2.3} and Proposition~\hl{P3.4}{3.4} imply that
\htt{3.7.2}{}
$${}^P\!n_{f,k/3m}=\mu^{(2)}_k=\mu'_k+\mu''_k-\nu_{k+3m}=\begin{cases}0&(\,k\les 1\,),\\
k-1&(\,1\les k\les m+1\,),\\
m&(\,m+1\les k\les 3m-1\,),\\
4m+1-k\h{ }&(\,3m\les k\les 4m+1\,),\\
0&(\,4m+1\les k\,).\end{cases}
\leqno(3.7.2)$$
\par\nin (Note that $\mu^{(2)}_k+\mu^{(2)}_{k+3m}=m$ ($k\in[1,3m-1])$ and $\mu^{(2)}_{3m}=m+1$.)
Indeed, we have by Lemma~\hl{L2.3}{2.3}
$$\aligned\mu'_k&=\begin{cases}0&(\,k\les 2m-1\,),\\
k-2m+1\h{ }&(\,2m-1\les k\les 3m\,),\\
4m+1-k&(\,3m\les k\les 4m+1\,),\\
0&(\,4m+1\les k\,),\end{cases}\\
\mu''_k&=\begin{cases}0&(\,k\les 1\,),\\
k-1&(\,1\les k\les 2m-1\,),\\
2m-2\q\q\!\h{ }&(\,2m-1\les k\,),\end{cases}\\
\nu_{k+3m}&=\begin{cases}0&(\,k\les m+1\,),\\
k-m-1\,\,\,\h{ }&(\,m+1\les k\les 3m-1\,),\\
2m-2&(\,3m-1\les k\,).\end{cases}\endaligned$$
\par\nin These formulas show that the relation between the Steenbrink spectrum $\Sp(f)$ and the pole order spectrum $\Sp_P(f)$ is rather complicated even for $n=2$ in general.
\bs\bs\htt{S4}{}
\vbox{\centerline{\bf 4. Gauss-Manin systems and Brieskorn modules}
\bsn
In this section we prove Theorem~\hl{T2}{2} after recalling some facts from Gauss-Manin systems and Brieskorn modules.}
\msn
{\bf 4.1.~Graded Gauss-Manin complexes.} Let $f$ be a homogeneous polynomial in $R$ with degree $d$. In the notation of \hl{1.1}{1.1}, the graded Gauss-Manin complex $C_f^{\ssb}$ associated with $f$ is defined by
$$C_f^j:=\Om^j[\dd_t]\q(j\in\Z),$$
\par\nin where $\dd_t$ has degree $-d$. This means that
$$\Om^j\,\dd_t^p=\Om^j(pd),$$
\par\nin where $(pd)$ denotes the shift of the grading as in the introduction. Its differential $\ddd$ is defined by
\htt{4.1.1}{}
$$\ddd(\omega\,\dd_t^p)=(\ddd\omega)\,\dd_t^p-(\df\wedge\omega)\,
\dd_t^{p+1}\q\h{for}\,\,\,\omega\in\Om^k.
\leqno(4.1.1)$$
\par\nin where $\ddd\omega$ denotes the differential of the de Rham complex. It has a structure of a complex of $\C[t]\langle\dd_t\rangle$-modules defined by
\htt{4.1.2}{}
$$t(\omega\,\dd_t^p)=(f\omega)\,\dd_t^p-p\,\omega\,\dd_t^{p-1},\q\dd_t (\omega\,\dd_t^p)=\omega\,\dd_t^{p+1}\q\h{for}\,\,\,\omega\in\Om^j.
\leqno(4.1.2)$$
\par\nin The Gauss-Manin systems are defined by the cohomology groups $H^jC_f^{\ssb}\,\,(j\in\Z)$. These are regular holonomic graded $\C[t]\langle\dd_t\rangle$-modules. By the same argument as in \cite{BaSa}, we have
\htt{4.1.3}{}
$$\h{The action of $\dd_t$ on $H^jC_f^{\ssb}$ is bijective for $j\ne 1$.}
\leqno(4.1.3)$$
\par\nin \par\htt{4.2}{}\msn
{\bf 4.2.~Brieskorn modules.} Let $(A_f^{\ssb},\ddd)$ be a graded subcomplex of the de Rham complex $(\Om^{\ssb},\ddd)$ defined by
$$A_f^j:={\rm Ker}(\df\wedge:\Om^j\to\Om^{j+1}(d)).$$
\par\nin The Brieskorn modules are graded $\C[t]\langle\dd_tt\rangle$-modules defined by its cohomology groups
$$H^jA_f^{\ssb}\q(j\in\Z).$$
\par\nin The actions of $t$, $\dti$, $\dd_tt$ are respectively defined by the multiplication by $f$,
$$\aligned\dti[\omega]=[\df\wedge\xi]\q\h{with}\q\ddd\xi=\omega,\\
\dd_tt\,[\omega]=[\ddd\eta]\q\h{with}\q\df\wedge\eta=f\omega,\endaligned$$
\par\nin where $[\omega]$ denotes the cohomology class, see \cite{Bri}, \cite{BaSa}, etc. (In case $j=1$, we have to choose a good $\xi$ for the action of $\dti$, see \cite{BaSa}.) Moreover, we have
\htt{4.2.1}{}
$$\dd_tt\,[\omega]=(k/d)[\omega]\q\h{for}\q[\omega]\in(H^jA_f^{\ssb})_k,
\leqno(4.2.1)$$
\par\nin where $(H^jA_f^{\ssb})_k$ denotes the degree $k$ part. (This follows from the definition by using the contraction with the Euler vector field $\xi:=\msum_i\,x_i\,\dd/\dd x_i$.) This implies
\htt{4.2.2}{}
$$t\,[\omega]=(k/d)\,\dti[\omega]\q\h{for}\q[\omega]\in(H^jA_f^{\ssb})_k.
\leqno(4.2.2)$$
\par\nin Since $(H^jA_f^{\ssb})_k=0$ for $k\les 0$, this implies that ${\rm Coker}^n\,t$ in Theorem~\hl{T2}{2} can be replaced with ${\rm Coker}^n\,\dti$.
\sk
There is a natural inclusion
$$A_f^{\ssb}\into C_f^{\ssb}.$$
\par\nin This is compatible with the actions of $t$, $\dti$, $\dd_tt$ on the cohomology by definition. So (\hl{4.2.1}{4.2.1}) holds also for $\omega\in(H^jC_f^{\ssb})_j$, since the image of $H^jA_f^{\ssb}$ generates $H^jC_f^{\ssb}$ over $\C[\dd_t]$. The last assertion is well known in the analytic case (see for instance \cite{BaSa}), and is reduced to this case by using the scalar extensions
$$R\into\C\{x_1,\dots,x_n\},\q
\C[t]\into\C\{t\}.$$
\par\nin \sk
For $j\in\Z$, we then get
\htt{4.2.3}{}
$$H^{j+1}(C_f^{\ssb})_k=\begin{cases}H^j(f^{-1}(1),\C)_{\la}&
\h{if}\,\,\,k/d\notin\Z_{\les 0},\\
\HH^j(f^{-1}(1),\C)_{\la}&
\h{if}\,\,\,k/d\notin\Z_{\ges 1},\end{cases}
\leqno(4.2.3)$$
\par\nin in the notation of \hl{3.2}{3.2}, where $\la=\exp(-2\pi ik/d)$, see also \cite{Di1}.
\sk
We have moreover
\htt{4.2.4}{}
$${\rm Ker}(H^jA_f^{\ssb}\to H^jC_f^{\ssb})=(H^jA_f^{\ssb})_{\rm tor},
\leqno(4.2.4)$$
\par\nin where the last term denotes the $t$-torsion subspace of $H^jA_f^{\ssb}$, which coincides with the $\dti$-torsion, and is annihilated by $\dd_t^{-p}$ for $p\gg 0$, see \cite{BaSa}.
\msn
{\bf 4.3.~Relation with the Koszul cohomologies.} Set
\htt{4.3.1}{}
$$A_f^{\prime\,j}:=\df\wedge\Om^{j-1}\buildrel\iota\over\into A_f^j\q(j\in\Z).
\leqno(4.3.1)$$
\par\nin Using the short exact sequence of complexes
$$0\to(A_f^{\ssb},\ddd)\to(\Om^{\ssb},\ddd)\to (A_f^{\prime\,\ssb},\ddd)[1]\to 0,$$
\par\nin we get isomorphisms
\htt{4.3.2}{}
$$\dd:H^jA_f^{\prime\,\ssb}\simto H^jA_f^{\ssb}\q(j\ne 1),
\leqno(4.3.2)$$
\par\nin together with a short exact sequence
$$0\to\C\to H^1A_f^{\prime\,\ssb}\to H^1A_f^{\ssb}\to 0.$$
\par\nin By (\hl{4.3.1}{4.3.1}) and (\hl{4.3.2}{4.3.2}), we get an action of $\dti$ on $H^jA_f^{\prime\,\ssb}$, $H^jA_f^{\ssb}$ defined respectively by
$$\dti:=\dd^{-1}\ssc H^j\iota,\q\dti:=H^j\iota\ssc\dd^{-1}
\q(j\ne 1).$$
\par\nin \sk
We have the canonical isomorphism
\htt{4.3.3}{}
$$(A_f^{\ssb}/A_f^{\prime\,\ssb},\ddd)=(H^{\ssb}_{\ddd f\wedge}\Om^{\ssb},\ddd),
\leqno(4.3.3)$$
\par\nin where $H^{\ssb}_{\ddd f\wedge}$ means that the cohomology is taken for the differential $\ddd f\wedge$ of $\Om^{\ssb}$ (preserving the grading up to the shift by $-d$), and $\ddd$ acts on $H^{\ssb}_{\ddd f\wedge}$ by the anti-commutativity of $\ddd$ and $\df\wedge$.
The relation with the shifted Koszul complex $({}^s\!K_f^{\ssb},\df\wedge)$ in the introduction is given by 
$$\Hdf^{j+n}\Om^{\ssb}=H^j({}^s\!K_f^{\ssb})(-jd)\q(j\in[-n,0]).$$
\par\nin By the short exact sequence of complexes
$$0\to(A_f^{\prime\,\ssb},\ddd)\buildrel\iota\over\to(A_f^{\ssb},\ddd)\to(\Hdf^{\ssb}\Om^{\ssb},\ddd)\to 0,$$
\par\nin we get a long exact sequence
\htt{4.3.4}{}
$$\to H^{j-1}_{\ddd}(\Hdf^{\ssb}\Om^{\ssb})\to H^jA_f^{\prime\,\ssb}
\buildrel{\iota_j}\over\to H^jA_f^{\ssb}\to H^j_{\ddd}(\Hdf^{\ssb}\Om^{\ssb})\to,
\leqno(4.3.4)$$
\par\nin where $H^j_{\ddd}$ means that the cohomology is taken for the differential $\ddd$, and the middle morphism $\iota_j:=H^j\iota$ can be identified by (\hl{4.3.2}{4.3.2}) with
$$\dti:H^jA_f^{\ssb}\to H^jA_f^{\ssb}\,\,\,\h{if}\,\,\,j>1.$$
\par\nin In particular we get for $j=n$
\htt{4.3.5}{}
$$H^n_{\ddd}(\Hdf^{\ssb}\Om^{\ssb})={\rm Coker}(\dti:H^nA_f^{\ssb}\to H^nA_f^{\ssb}).
\leqno(4.3.5)$$
\par\nin \sk
By the above argument, the $\dti$-torsion of $H^jA_f^{\ssb}$ contributes to $H^{j-1}_{\ddd}(\Hdf^{\ssb}\Om^{\ssb})$, and we get in particular
\htt{4.3.6}{}
$$\h{$H^cA_f^{\ssb}$ is torsion-free if $c$ is the codimension of $\Sg\,f^{-1}(0)\subset\C^n$.}
\leqno(4.3.6)$$
\par\nin Note that $c=n-1$ under the assumption of the introduction. By Theorems~\hl{T5.2}{5.2} and \hl{T5.3}{5.3} below, the $\dti$-torsion of $H^nA_f^{\ssb}$ is finite dimensional if and only if all the singularities of $Z$ are weighted homogeneous.
\msn
{\bf 4.4.~Filtrations $P'$ and $G$.} There are two filtrations $P'$, $G$ on $C_f^{\ssb}$ defined by
\htt{4.4.1}{}
$$\aligned P'_p\,C_f^k&:=\mopl_{i\les k+p}\,\Om^k\,\dd_t^i,\\ G_p\,C_f^k&:=\bl(\mopl_{i<p}\,\Om^k\,\dd_t^i\br)\oplus A_f^p\,\dd_t^p.
\endaligned
\leqno(4.4.1)$$
\par\nin These are exhaustive increasing filtrations. Set $P^{\prime\,p}=P'_{-p}$, $G^p=G_{-p}$. By definition, we have
\htt{4.4.2}{}
$$\Gr^p_{P'}C_f^{\ssb}=\sigma_{\ges p}\bl(K_f^{\ssb}((n-p)d)\br),
\leqno(4.4.2)$$
\par\nin see \cite{De} for the truncation $\sigma_{\ges p}$. Let ${\rm Dec}\,P'$ be as in loc.~cit. Then we have
\htt{4.4.3}{}
$$G={\rm Dec}\,P'.
\leqno(4.4.3)$$
\par\nin Since the differential of $C_f^{\ssb}$ respect the grading, we have the pole order spectral sequence in the category of graded $\C$-vector spaces
\htt{4.4.4}{}
$$_{P'}E_1^{p,j-p}=H^j\Gr^p_{P'}C_f^{\ssb}\Longrightarrow H^jC_f^{\ssb},
\leqno(4.4.4)$$
\par\nin with
\htt{4.4.5}{}
$$_{P'}E_1^{p,j-p}=\begin{cases}\,0&\h{if}\,\,\,j<p,\\ A_f^p&\h{if}\,\,\,j=p,\\ H^jK_f^{\ssb}((n-p)d)&\h{if}\,\,\,j>p,\end{cases}
\leqno(4.4.5)$$
\par\nin \htt{4.4.6}{}
$$_{P'}E_2^{p,j-p}=\begin{cases}\,0&\h{if}\,\,\,j<p,\\ H^pA_f^{\ssb}&\h{if}\,\,\,j=p,\\ H^j_{\ddd}(\Hdf^{\ssb}\Om^{\ssb})((j-p)d)&\h{if}\,\,\,j>p,\end{cases}
\leqno(4.4.6)$$
\par\nin where $\Hdf^{\ssb}\Om^{\ssb}$ is as in (\hl{4.3.3}{4.3.3}).
\sk
Note that the degeneration at $E_2$ of the pole order spectral sequence is equivalent to the strictness of ${\rm Dec}\,P'$ by \cite{De}, and the latter condition is equivalent to the torsion-freeness of the $H^jA_f^{\ssb}$ by using (\hl{4.2.4}{4.2.4}) and (\hl{4.4.3}{4.4.3}). The obtained equivalence seems to be known to the specialists (see for instance \cite{vSt}), and the above argument may simplify some argument in loc.~cit.
\sk
By the isomorphism (\hl{4.2.3}{4.2.3}) for $k\in[1,d]$, the filtration $P'$ on the left-hand side of (\hl{4.2.3}{4.2.3}) induces a filtration $P'$ on the right-hand side. This corresponds to the filtration $P$ by the isomorphism (\hl{3.1.2}{3.1.2}) up to the shift of the filtration by 1, and we get the isomorphisms
\htt{4.4.7}{}
$$P^{\prime\,p+1}H^{j+1}(C_f^{\ssb})_k\cong P^pH^j(f^{-1}(1),\C)_{\la}\q\bl(\la=\exp(-2\pi ik/d),\,k\in[1,d]\br),
\leqno(4.4.7)$$
\par\nin see \cite[Ch.~6, Thm.~2.9]{Di1} (and also \cite[Section 1.8]{DiSa2} in case $j=n-1$). By (\hl{3.1.3}{3.1.3}), we have the inclusions
\htt{4.4.8}{}
$$F^p\subset P^p\q\h{on}\,\,\,\,H^j(f^{-1}(1),\C)_{\la},
\leqno(4.4.8)$$
\par\nin Here it is possible to show (\hl{4.4.8}{4.4.8}) by calculating the direct image of $(\OO_X,F)$ by $f$ as a filtered $\Dc$-module underlying a mixed Hodge module, see \cite{mhp}, \cite{mhm}, where a compactification of $f$ must be used. (The shift of the filtration by 1 comes from the direct image of $\OO_X$ as a {\it left} $\Dc$-module by the graph embedding of $f$.)
\sk
The inclusion (\hl{4.4.8}{4.4.8}) implies some relation between the spectrum and the Poincar\'e series of the Koszul cohomologies via the spectral sequence (\hl{4.4.4}{4.4.4}), and the difference between $F^p$ and $P^p$ implies also their difference in certain cases, see also \cite{Di1}, \cite{Di3}, \cite{DiSt2}.
\par\htt{4.5}{}\msn
{\bf 4.5.~Algebraic microlocal Gauss-Manin complexes.} For a homogeneous polynomial $f$, let $\Ct^{\ssb}_f$ be the algebraic microlocal Gauss-Manin complex (that is, $\Ct_f^j=\Om^j[\dd_t,\dti]$). The algebraic microlocal Gauss-Manin systems $H^j\Ct^{\ssb}_f$ are free graded $\C[\dd_t,\dti]$-modules of finite type. Replacing $C^{\ssb}_f$ with $\Ct^{\ssb}_f$ in (\hl{4.4.1}{4.4.1}) and (\hl{4.4.4}{4.4.4}), we have the filtrations $P'$, $G$ on $\Ct^{\ssb}_f$ together with the microlocal pole order spectral sequence
\htt{4.5.1}{}
$$_{P'}\E_1^{p,j-p}=H^j\Gr^p_{P'}\Ct_f^{\ssb}\Longrightarrow H^j\Ct_f^{\ssb},
\leqno(4.5.1)$$
\par\nin where (\hl{4.4.3}{4.4.3}) holds again (that is, $G={\rm Dec}\,P'$), and the last equalities of (\hl{4.4.5}{4.4.5}) and (\hl{4.4.6}{4.4.6}) hold for any $j,p\in\Z$, that is,
\htt{4.5.2}{}
$${}_{P'}\E_r^{p,j-p}=\begin{cases}\Hdf^j\Om^{\ssb}((j-p)d)=H^jK_f^{\ssb}((n-p)d)&\h{if}\,\,\,r=1,\\
H^j_{\ddd}(\Hdf^{\ssb}\Om^{\ssb})((j-p)d)&\h{if}\,\,\,r=2.\end{cases}
\leqno(4.5.2)$$
\par\nin Moreover the last equality of (\hl{4.2.3}{4.2.3}) holds for any $k$, that is,
\htt{4.5.3}{}
$$H^{j+1}(\Ct_f^{\ssb})_k=\HH^j(f^{-1}(1),\C)_{\la}\q\h{with}\q\la=\exp(-2\pi ik/d),
\leqno(4.5.3)$$
\par\nin (Note that the Gauss-Manin complex $C_f^{\ssb}$ can be defined also as the single complex associated with the double complex having two differentials $\ddd$ and $\df\wedge$, see \cite{Di1}, \cite{Di3}, etc.)
\sk
Let $P',G$ denote also the induced filtrations on $H^j(C_f^{\ssb})$, $H^j(\Ct_f^{\ssb})$ .There is a canonical inclusion
$$C_f^{\ssb}\into\Ct_f^{\ssb}.$$
\par\nin Set
$$\omega_0:=\df\in H^1(G_0C_f^{\ssb})\,(=H^1A_f^{\ssb}).$$
\par\nin By the same argument as in \cite{BaSa}, it generates a free $\C[t]$-module for $p\in\N\cup\{\infty\}$
$$\C[t]\omega_0\subset H^1(G_pC_f^{\ssb}),$$
\par\nin where $G_{\infty}C_f^{\ssb}:=C_f^{\ssb}$. 
Set
$$\HH^j(G_pC_f^{\ssb})=\begin{cases}H^j(G_pC_f^{\ssb})&\h{if}\,\,\,j\ne 1,\\ H^j(G_pC_f^{\ssb})/\C[t]\omega_0&\h{if}\,\,\,j=1.\end{cases}$$
\par\nin Then the above inclusion induces the canonical isomorphisms
\htt{4.5.4}{}
$$\HH^j(G_pC_f^{\ssb})\simto H^j(G_p\Ct_f^{\ssb})\q(p\in\N\cup\{\infty\},\,\,j\in\Z).
\leqno(4.5.4)$$
\par\nin Indeed, the assertion for $p=\infty$ follows from the same argument as in loc.~cit. This implies the assertion for $p\in\N$ by using the canonical morphism of long exact sequences
$$\begin{CD}@>>>\HH^j(G_pC_f^{\ssb})@>>>\HH^j(C_f^{\ssb})@>>>
H^j(C_f^{\ssb}/G_pC_f^{\ssb})@>>>\\
@. @VVV @VVV @|\\
@>>>H^j(G_p\Ct_f^{\ssb})@>>>H^j(\Ct_f^{\ssb})@>>>
H^j(\Ct_f^{\ssb}/G_p\Ct_f^{\ssb})@>>>\end{CD}$$
\par\nin \sk
From the canonical isomorphisms (\hl{4.5.4}{4.5.4}), we can deduce
\htt{4.5.5}{}
$$G_p\HH^j(C_f^{\ssb})\simto G_pH^j(\Ct_f^{\ssb})=\dd_t^pG_0H^j(\Ct_f^{\ssb})\q(p\in\N,\,\,j\in\Z).
\leqno(4.5.5)$$
\par\nin This implies
\htt{4.5.6}{}
$$\dd_t:\Gr^G_p\HH^j(C_f^{\ssb})_k\simto\Gr^G_{p+1}\HH^j(C_f^{\ssb})_{k-d}\q(p\in\N,\,\,j,k\in\Z).
\leqno(4.5.6)$$
\par\nin Note that these hold with $G$ replaced by $P'$ by (\hl{4.4.3}{4.4.3}).
We then get by (\hl{4.4.7}{4.4.7})
\htt{4.5.7}{}
$$P^{\prime\,p+1}H^{j+1}(\Ct_f^{\ssb})_k\cong P^p\HH^j(f^{-1}(1),\C)_{\la}\q\bl(\la=\exp(-2\pi ik/d),\,k\in[1,d]\br),
\leqno(4.5.7)$$
\par\nin \par\htt{P4.6}{}\msn
{\bf Proposition~4.6.} {\it With the notation of $(4.4)$ and $(4.5)$, there are canonical isomorphisms for $r\ges 2$
$$\aligned{\rm Im}(\ddd_r:{}_{P'}E_r^{\,p-r,n-p+r-1}\to{}_{P'}E_r^{\,p,n-p})&=\begin{cases}\,0&\h{if}\,\,\,p>n,\\\Gr^K_{r-1}(H^nA_f^{\ssb})_{\rm tor}&\h{if}\,\,\,p=n,\\ \Gr^K_{r-1}({\rm Coker}^n\,\dti)((n-p)d)&\h{if}\,\,\,p<n,\end{cases}\\
{\rm Im}(\ddd_r:{}_{P'}\E_r^{\,p-r,n-p+r-1}\to{}_{P'}\E_r^{\,p,n-p})&=\Gr^K_{r-1}({\rm Coker}^n\,\dti)((n-p)d),\endaligned$$
\par\nin where $K_{\ssb}$ is the kernel filtration, and ${\rm Coker}^n\,\dti$ is a quotient of $(H^nA_f^{\ssb})_{\rm tor}$ as in Theorem~$\hl{T2}{2}$.}
\msn
{\it Proof.} We first show the assertion for the microlocal pole order spectral sequence, that is, for the second isomorphism. Since ${}_{P'}\E_r^{p,j-p}=0$ for $j>n$, the images of the differentials
$$\ddd_r:{}_{P'}\E_r^{\,p-r,n-p+r-1}\to{}_{P'}\E_r^{\,p,n-p}\,\,\,(r\ges 2)$$
\par\nin correspond to an increasing sequence of subspaces (with $p$ fixed):
\htt{4.6.1}{}
$$\I_r^{\,p,n-p}\subset{}_{P'}\E_2^{\,p,n-p}=({\rm Coker}^n\,\dti)((n-p)d)\,\,\,(r\ges 2),
\leqno(4.6.1)$$
\par\nin such that
$${\rm Im}\bl(\ddd_r:{}_{P'}\E_r^{\,p-r,n-p+r-1}\to{}_{P'}\E_r^{\,p,n-p}\br)=\I_r^{\,p,n-p}/\I_{r-1}^{\,p,n-p}\,\,\,(r\ges 2),$$
\par\nin with $\I_1^{\,p,n-p}:=0$. Here ${\rm Coker}^n\,\dti$ is a quotient of $H^nA_f^{\ssb}$ (and not $(H^nA_f^{\ssb})_{\rm tor})$, and (\hl{4.3.5}{4.3.5}) is used for the last isomorphism of (\hl{4.6.1}{4.6.1}).
\sk
By the construction of the spectral sequence (see for instance \cite{De}), we have
\htt{4.6.2}{}
$$\I_r^{\,p,n-p}=K_{r-1}({\rm Coker}^n\,\dti)((n-p)d),
\leqno(4.6.2)$$
\par\nin where $K_{\ssb}$ is the kernel filtration defined just before Theorem~\hl{T2}{2}.
(More precisely, $K_{\ssb}$ defines a non-exhaustive filtration of $H^nA_f^{\ssb}$, and its union is $(H^nA_f^{\ssb})_{\rm tor}$.)
Indeed, the left-hand side is given by the classes of $\omega\in\Om^n$ such that there are
$$\eta_i\in\Om^n\,\,(i\in[0,r-1])$$
\par\nin satisfying
$$d\eta_0=\omega,\q d\eta_{i+1}=df\wedge\eta_i\,(i\in[0,r-2]),\q df\wedge\eta_{r-1}=0.$$
\par\nin However, this condition is equivalent to that the class of $\omega$ in the Brieskorn module is contained in $K_{r-1}(H^nA_f^{\ssb})_{\rm tor}$.
(Note that $[df\wedge\eta_{r-2}]$ gives $\dd_t^{1-r}[\omega]$ and vanishes in $H^nA_f^{\ssb}$.)
So the second isomorphism follows.
\sk
The argument is essentially the same for the first isomorphism by replacing (\hl{4.6.2}{4.6.2}) with
$$I_r^{\,p,n-p}=\begin{cases}\,0&\h{if}\,\,\,p>n,\\ K_{r-1}(H^nA_f^{\ssb})_{\rm tor}&\h{if}\,\,\,p=n,\\ K_{r-1}({\rm Coker}^n\,\dti)((n-p)d)&\h{if}\,\,\,p<n.\end{cases}$$
\par\nin This finishes the proof of Proposition~\hl{P4.6}{4.6}.
\ms
As a corollary of Proposition~\hl{P4.6}{4.6}, we get the following.
\par\htt{C4.7}{}\msn
{\bf Corollary~4.7.} {\it The following three conditions are equivalent to each other$\,:$
\skn
$(a)$ The pole order spectral sequence $(4.4.4)$ degenerates at $E_2$.
\skn
$(b)$ The algebraic microlocal pole order spectral sequence $(4.5.1)$ degenerates at $E_2$.
\skn
$(c)$ The torsion subgroup $(H^nA_f^{\ssb})_{\rm tor}$ vanishes.}
\msn
{\bf 4.8.~Proof of Theorem~\hl{T2}{2}.} By (\hl{4.5.7}{4.5.7}) the assertion follows from the second isomorphism in Proposition~\hl{P4.6}{4.6} by choosing any $p\in\Z$, where the obtained isomorphism is independent of the choice of $p$ by using the bijectivity of the action of $\dd_t$. (It is also possible to use the first isomorphism in Proposition~\hl{P4.6}{4.6} by choosing some $p<n$ although the independence of the choice of $p$ is less obvious unless the relation with the algebraic microlocal pole order spectral sequence is used.) This finishes the proof of Theorem~\hl{T2}{2}.
\par\htt{4.9}{}\msn
{\bf 4.9.~Thom-Sebastiani type theorem for $P'$.} Let $f,f_1,f_2$ be as in Proposition~\hl{P2.2}{2.2}. In the notation of \hl{4.5}{4.5}, we have a canonical isomorphism
$$(\Ct^{\ssb}_f,P')=(\Ct^{\ssb}_{f_1},P')\otimes_{\C[\dd_t,\,\dd_t^{-1}]}(\Ct^{\ssb}_{f_2},P').$$
\par\nin Assume $f_2$ has an isolated singularity at the origin as in Proposition~\hl{P2.2}{2.2}. Then
$$H^jGr^{P'}_k\Ct^{\ssb}_{f_2}=0\q(j\ne n_2,\,\,k\in\Z).$$
\par\nin Hence $(\Ct^{\ssb}_{f_2},P')$ is strict, and we get a filtered quasi-isomorphism
$$(\Ct^{\ssb}_{f_2},P')\simto H^{n_2}(\Ct^{\ssb}_{f_2},P')[-n_2].$$
\par\nin This implies a filtered quasi-isomorphism
\htt{4.9.1}{}
$$(\Ct^{\ssb}_f,P')\simto(\Ct^{\ssb}_{f_1},P')\otimes_{\C[\dd_t,\,\dd_t^{-1}]} H^{n_2}(\Ct^{\ssb}_{f_2},P')[-n_2],
\leqno(4.9.1)$$
\par\nin which is compatible with the action of $t$. More precisely, the action of $t$ on the left-hand side corresponds to $t\otimes id+id\otimes t$ on the right-hand side (since $f=f_1+f_2$).
\sk
Combining (\hl{4.9.1}{4.9.1}) with (\hl{4.5.7}{4.5.7}), we get the Thom-Sebastiani type theorem for the pole order spectrum:
\htt{4.9.2}{}
$$\Sp_P(f)=\Sp_P(f_1)\,\Sp_P(f_2),\q\Sp_P^j(f)=\Sp_P^j(f_1)\,\Sp_P^0(f_2)\q(j\in\N),
\leqno(4.9.2)$$
\par\nin assuming that $f_2$ has an isolated singularity as above so that $\Sp_P(f_2)=\Sp_P^0(f_2)$, see \cite{ScSt} for the case where $f_1$ has also an isolated singularity.
Note that the Thom-Sebastiani type theorem holds for the Steenbrink spectrum by \cite{ts}.
\msn
{\bf Remarks~4.10.} (i) With the notation and assumption of \hl{4.9}{4.9}, the pole order spectral sequences degenerate at $E_2$ for $f$ if and only if they do for $f_1$. This follows from (\hl{4.9.1}{4.9.1}) together with Corollary~\hl{C4.7}{4.7}.
\ms
(ii) The equivalence between the $E_2$-degeneration of the pole order spectral sequence (\hl{4.4.4}{4.4.4}) and the vanishing of $(H^nA_f^{\ssb})_{\rm tor}$ was shown in \cite{vSt} in the (non-graded) analytic local case.
\ms
(iii) Assuming $\dim\Sg\,f^{-1}(0)=1$, we have by (\hl{4.3.4}{4.3.4}) the following exact sequence:
$$\aligned
0\to\HH^{n-1}A_f^{\ssb}(-d)&\,{\buildrel{\dti}\over\longrightarrow}\,\HH^{n-1}A_f^{\ssb}\to H^{n-1}_{\ddd}(\Hdf^{\ssb}\Om^{\ssb})\\
\to H^nA_f^{\ssb}(-d)&\,{\buildrel{\dti}\over\longrightarrow}\,H^nA_f^{\ssb}\to H^n_{\ddd}(\Hdf^{\ssb}\Om^{\ssb})\to 0,\endaligned$$
\par\nin where $\HH^{n-1}A_f^{\ssb}$ is defined by $H^{n-1}A_f^{\ssb}$ if $n\ne2$, and by its quotient by $\C[t]\omega_0$ if $n=2$. (For $\omega_0$, see the definition of $\HH^j(G_pC_f^{\ssb})$ in \hl{4.5}{4.5}.)
This exact sequence has sufficient information about the torsion subgroup $H^n_{\ddd}(\Hdf^{\ssb}\Om^{\ssb})_{\rm tor}$ to give another proof of Theorem~\hl{T2}{2}.
\ms
(iv) By forgetting the grading, Proposition~\hl{P4.6}{4.6} and Corollary~\hl{C4.7}{4.7} can be extended to the analytic local case where $f$ is a germ of a holomorphic function on a complex manifold with $\dim\Sg\,f^{-1}(0)=1$.
\ms
The following will be used in the proof of Theorem~\hl{T5.2}{5.2} below.
\par\htt{4.11}{}\msn
{\bf 4.11.~Multiplicity of the minimal exponent.}
Let $g$ be a germ of holomorphic function on a complex manifold $(Y,0)$ having an isolated singularity. We have the direct image $\B:=\OO_{Y,0}[\dd_t]$ of $\OO_{Y,0}$ as a left $\Dc_{Y,0}$-module by the graph embedding of $g$. (Note that it is an analytic $\Dc$-module.) It has the Hodge filtration $F$ by the order of $\dd_t$ and the filtration $V$ of Kashiwara \cite{Ka} and Malgrange \cite{Ma}.
\sk
Consider $\Gr_V^{\alpha}(\B,F)$ for $\alpha<1$. These underlie mixed Hodge modules supported at $0$, and are the direct images of filtered vector spaces by the inclusion $\{0\}\into Y$ as filtered $\Dc$-modules. (This is shown by using \cite[Lemma 3.2.6]{mhp} applied to any function vanishing at $0$.) So we get
\htt{4.11.1}{}
$$\h{The $\Gr^F_p\Gr_V^{\alpha}\B$ are annihilated by ${\mathfrak m}_{Y,0}\subset\OO_{Y,0}$ for $\alpha<1$,}
\leqno(4.11.1)$$
\par\nin where ${\mathfrak m}_{Y,0}\subset\OO_{Y,0}$ is the maximal ideal.
\sk
Let $\BB:=\OO_{Y,0}[\dd_t,\dti]$ be the algebraic microlocalization of $\B$. By \cite[Sections 2.1-2]{mhm}, it has the Hodge filtration $F$ by the order of $\dd_t$ and also the filtration $V$ such that
$$\aligned\dd_t:F_pV^{\alpha}\BB&\simto F_{p+1}V^{\alpha-1}\BB\q(\forall\,p,\alpha).\\(\Gr_V^{\alpha}\B,F)&\simto(\Gr_V^{\alpha}\BB,F)\q(\alpha<1),\endaligned$$
\par\nin Then (\hl{4.11.1}{4.11.1}) implies
\htt{4.11.2}{}
$$\h{The $\Gr^F_p\Gr_V^{\alpha}\BB$ are annihilated by ${\mathfrak m}_{Y,0}\subset\OO_{Y,0}$ for any $\alpha$.}
\leqno(4.11.2)$$
\par\nin \sk
Consider the (relative) de Rham complexes
$$\CC:={\rm DR}_Y(\B),\q\CCt:={\rm DR}_Y(\BB).$$
\par\nin Up to a shift of complexes, these are the Koszul complexes associated with the action of $\dd_{y_i}$ on $\B$ and $\BB$ where the $y_i$ are local coordinates of $Y$. It has the filtrations $F$ and $V$ induced by those on $\B$ and $\BB$. Here $V$ is stable by the action of $\dd_{y_i}$, but we need a shift for $F$ depending on the degree of the complexes $\CC$, $\CCt$. By the above argument we have
\htt{4.11.3}{}
$$H^j\Gr_p^F\Gr_V^{\alpha}\CCt=H^j\Gr_V^{\alpha}\CCt=H^j\Gr_p^F\CCt=0\q(j\ne 0),
\leqno(4.11.3)$$
\par\nin where we also use the assertion that $\Gr_p^F\CCt$ is the Koszul complex for the regular sequence $\{\dd g/\dd y_j\}$. These imply the vanishing of $H^jF_p\Gr_V^{\alpha}\CCt$, etc.\ for $j\ne 0$, and we get
\htt{4.11.4}{}
$$\h{$(\CCt;F,V)$ is strict,}
\leqno(4.11.4)$$
\par\nin by showing the exactness of some commutative diagram appearing in the definition of strict complex \cite{mhp}.
\sk
It is known that the filtration $V$ on $\CC$ is strict, and induces the filtration $V$ of Kashiwara and Malgrange on the Gauss-Manin system $H^0\CC$ (by using the arguments in the proof of \cite[Prop.~3.4.8]{mhp}).
This assertion holds by replacing $\CC$ with $\CCt$, since $\CC/V^{\alpha}\CC=\CCt/V^{\alpha}\CCt$ for $\alpha\les 1$ and $H^0\CC=H^0\CCt$ (see for instance \cite{BaSa}). Here we also get the canonical isomorphism
\htt{4.11.5}{}
$$(H^0\CC,V)=(H^0\CCt,V).
\leqno(4.11.5)$$
\par\nin \sk
Consider now $(\Gr^F_0\CCt,V)$. This is a complex of filtered $\OO_{Y,0}$-modules, and is strict. By the above argument we get the canonical isomorphism of filtered $\OO_{Y,0}$-modules
\htt{4.11.6}{}
$$H^0(\Gr^F_0\CCt,V)=(\OO_{Y,0}/(\dd g),V).
\leqno(4.11.6)$$
\par\nin Combining this with (\hl{4.11.2}{4.11.2}), (\hl{4.11.4}{4.11.4}) and using $\Gr_V^{\alpha}\Gr^F_p\CCt=\Gr^F_p\Gr_V^{\alpha}\CCt$, we get
\htt{4.11.7}{}
$$\h{The $\Gr_V^{\alpha}(\OO_{Y,0}/(\dd g))$ are annihilated by ${\mathfrak m}_{Y,0}\subset\OO_{Y,0}$ for any $\alpha$.}
\leqno(4.11.7)$$
\par\nin In particular, the multiplicity of the minimal exponent is 1.
\bs\bs\htt{S5}{}
\vbox{\centerline{\bf 5. Calculation of $\ddd^{(1)}$.}
\bsn
In this section we calculate $\ddd^{(1)}$ in certain cases, and prove Theorems~\hl{T5.2}{5.2} and \hl{T5.3}{5.3}.}
\par\htt{5.1}{}\msn
{\bf 5.1.~Relation with the isolated singularities in $\PP^{n-1}$.}
Let $\rho:\X\to X$ be the blow-up of the origin of $X:=\C^n$. Let $y=\sum_ic_ix_i$ be as in the introduction (that is, $(c_i)\in\C^n$ are sufficiently general). We may assume that
$$y=x_n,$$
\par\nin replacing the coordinates $x_1,\dots,x_n$ of $X=\C^n$. Let $\X'$ be the complement of the proper transform of $\{x_n=0\}$. It has the coordinates $\x_1,\dots,\x_n$ such that
$$\rho^*x_i=\begin{cases}\x_i\,\x_n&\h{if}\,\,\,i\ne n,\\
\x_n&\h{if}\,\,\,i=n.\end{cases}$$
\par\nin Define the complex $\spKf$ similarly to $\sKf$ in the introduction by replacing $R$ and $f$ respectively with
$$\C[\x_1,\dots,\x_{n-1}][\x_n,\x_n^{\,-1}]\q\h{and}\q f':=\rho^*f|_{\X'}=\x_n^d\,h(\x_1,\dots,\x_{n-1}).$$
\par\nin Here $h:=\rho^*f|_{\x_n=1}$. This is identified with $f|_{y=1}$, that is,
$$h(\x_1,\dots,\x_{n-1})=f(\x_1,\dots,\x_{n-1},1).$$
\par\nin Note that the grading of $\spKf$ is given only by the degree of $\x_n$ and $\ddd\x_n$.
\sk
The above construction of $\spKf$ is compatible with that of $\sKf$ via $\rho^*$, and we have the canonical graded morphism
$$H^j(\sKf)\to H^j(\spKf),$$
\par\nin in a compatible way with the differential $\ddd$. This morphism induces the injective morphisms
\htt{5.1.1}{}
$$N\into H^{-1}(\spKf),\q M''\into H^0(\spKf),
\leqno(5.1.1)$$
\par\nin where the image of $M'$ in $H^0(\spKf)$ vanishes. We have the inclusion
\htt{5.1.2}{}
$$N^{(2)}_{p+d}\subset{\rm Ker}\bl(\ddd:H^{-1}(\spKf)\to H^0(\spKf)\br)\cap N_{p+d},
\leqno(5.1.2)$$
\par\nin under the first injection of (\hl{5.1.1}{5.1.1}), and the equality holds if $M'_p=0$.
\sk
Let $Y'\,(\cong\C^{n-1})$ be the complement of $\{\x_n=0\}$ in $Y:=\PP^{n-1}$. Then
$$\X'=Y'\times\C,$$
\par\nin where $\x_1,\dots,\x_{n-1}$ and $\x_n$ are respectively coordinates of $Y'$ and $\C$. Moreover $\spKf$ is quasi-isomorphic to the mapping cone of
$$\dd f'/\dd \x_n=d\,\x_n^{\,d-1}h:(\Om_{Y'}^{n-1}/\ddh\wedge\Om_{Y'}^{n-2})[\x_n,\x_n^{\,-1}]\to (\Om_{Y'}^{n-1}/\ddh\wedge\Om_{Y'}^{n-2})[\x_n,\x_n^{\,-1}],$$
\par\nin where $\Om_{Y'}^j$ is identified with the group of global sections.
\sk
Let $\{z_i\}$ ($i\in I$) be the singular points of the morphism $h:Y'\to\C$. These are isolated singular points. (Indeed, they are the union of the singular points of $Y'_c:=\{h=c\}\subset Y'$ for $c\in\C$. Here the subvariety $Y'_c\subset Y'\,(=\C^{n-1})$ can be identified with the intersection of $\{f=c\}$ and $\{x_n=1\}$ in $X=\C^n$. Then the intersection of the closure of $Y'_c$ in $\PP^n\supset\C^n=X$ with the boundary $\PP^{n-1}=\PP^n\setminus\C^n$ is identified with the intersection of $Z=\{f=0\}$ and $\{x_n=0\}$ in $Y=\PP^n$, which is smooth by hypothesis. This implies that the singular points of $Y'_c$ form a proper variety contained in an affine variety. Hence they are discrete.)
\sk
Since the support of the $\C[\x_1,\dots,\x_{n-1}]$-module $\Om_{Y'}^{n-1}/\ddh\wedge\Om_{Y'}^{n-2}$ is $\{z_i\}_{i\in I}$, we have the canonical isomorphism
\htt{5.1.3}{}
$$\Om_{Y'}^{n-1}/\ddh\wedge\Om_{Y'}^{n-2}=\mopl_{i\in I}\,\Xi_{h_i}\q \h{with}\q\Xi_{h_i}:=\Om_{Y'_{\rm an},z_i}^{n-1}/\ddh_i\wedge\Om_{Y'_{\rm an},z_i}^{n-2},
\leqno(5.1.3)$$
\par\nin where $Y'_{\rm an}$ is the associated analytic space, and $(h_i,z_i):=(h_{\rm an},z_i)$.
\sk
Let $z_i$ ($i\in I_0$) be the singular points contained in $\{h=0\}$. These are the singular points of $Z:=\{f=0\}\subset\PP^{n-1}$, since $x_n$ is sufficiently general. For $i\notin I_0$, the analytic functions $h_i$ are invertible. Hence
\htt{5.1.4}{}
$$H^n(\spKf)=\mopl_{i\in I_0}\,(\Xi_{h_i}/h_i\,\Xi_{h_i})\wedge\C[\x_n,\x_n^{\,-1}]\,\ddd \x_n,
\leqno(5.1.4)$$
\par\nin and a similar formula holds for $H^{n-1}(\spKf)$ (with $\ddd \x_n$ on the right-hand side deleted, and $\wedge$ replaced by $\otimes$). 
So the $z_i$ ($i\notin I_0$) may be forgotten from now on.
\sk
By (\hl{5.1.4}{5.1.4}) and using the last inclusion of (\hl{5.1.1}{5.1.1}), we get
\htt{5.1.5}{}
$$M''\supset\mopl_{i\in I_0}\,(\Xi_{h_i}/h_i\,\Xi_{h_i})\wedge\C[\x_n]\,\x_n^{\,p}\,\ddd \x_n\q\h{for}\q p\gg 0.
\leqno(5.1.5)$$
\par\nin \sk
Take an element of pure degree $p$ of
$${\rm Ker}\bl(h_i:\Xi_{h_i}[\x_n,\x_n^{\,-1}]\to\Xi_{h_i}[\x_n,\x_n^{\,-1}]\br)\q(i\in I_0).$$
\par\nin It is represented by $\psi:=\frac{1}{d}\,\x_n^{\,p}\xi$ where $\xi\in\Om_{Y'_{\rm an},z_i}^{n-1}$ satisfies
\htt{5.1.6}{}
$$h_i\,\xi=\ddh_i\wedge\eta\q\h{with}\q\eta\in\Om_{Y'_{\rm an},z_i}^{n-2}.
\leqno(5.1.6)$$
\par\nin The corresponding element of $H^{n-1}(\spKf)$ is represented by
$$\psi':=\tfrac{1}{d}\,\x_n^{\,p}\xi+\x_n^{\,p-1}\ddd \x_n\wedge\eta.$$
\par\nin Its image in $H^n(\spKf)$ by the differential $\ddd$ is given by
\htt{5.1.7}{}
$$\ddd[\psi']=\pm\bl[\bl(\tfrac{p}{d}\,\xi-\ddd\eta\br)\wedge \x_n^{\,p-1}\ddd \x_n\br],
\leqno(5.1.7)$$
\par\nin and we have by (\hl{5.1.6}{5.1.6})
\htt{5.1.8}{}
$$[d\eta]=\dd_tt\,[\xi]\q\h{in}\,\,\,H''_{h_i}.
\leqno(5.1.8)$$
\par\nin \sk
Let $V$ be the $V$-filtration of Kashiwara \cite{Ka} and Malgrange \cite{Ma} on the Gauss-Manin system $\G_{h_i}$ indexed by $\Q$, see for instance \cite{ScSt}. (It is closely related to the theory of asymptotic Hodge structure \cite{Va}.) We denote also by $V$ the induced filtration on the Brieskorn module $H''_{h_i}$ and also on $\Xi_{h_i}$ via the canonical inclusion and the surjection
$$\G_{h_i}\supset H''_{h_i}\onto\Xi_{h_i},$$
\par\nin see \cite{Bri} for the latter. In this paper we index $V$ so that $\dd_tt-\alpha$ is nilpotent on $\Gr_V^{\alpha}\G_{h_i}$.
\sk
Let $\{\alpha_{h_i,l}\}$ be the exponents of $h_i$ counted with multiplicity; more precisely
\htt{5.1.9}{}
$$\#\{l:\alpha_{h_i,l}=\alpha\}=\dim\Gr^{\alpha}_V\Xi_{h_i}\q\h{and}\q\Sp_{h_i}(t)=\msum_l\,t^{\,\alpha_{h_i,l}}.
\leqno(5.1.9)$$
\par\nin Here we may assume the $\alpha_{h_i,l}$ are weakly increasing (that is, $\alpha_{h_i,l}\les\alpha_{h_i,l+1}$) for each $l$. We have the symmetry $\{\alpha_{h_i,l}\}_l=\{n-\alpha_{h_i,l}\}_l$ (counted with multiplicity) by \cite{St2}.
\par\htt{T5.2}{}\msn
{\bf Theorem~5.2.} {\it With the notation of {\rm\hl{5.1}{5.1}}, assume $h_i$ is non-quasihomogeneous $($that is, $h_i\notin(\dd h_i))$ for some $i\in I_0$. Then the kernel and cokernel of $d^{(1)}:N\to M$ $($that is, $N^{(2)}$ and $M^{(2)}$ in Theorem~$\hl{T2}{2})$ and $(H^nA_f^{\ssb})_{\rm tor}$ are all infinite dimensional over $\C$.}
\msn
{\it Proof.} Since the minimal exponent $\alpha_{h_i,1}$ has multiplicity 1 (see \hl{4.11}{4.11}), we have
$$V^{>\alpha_{h_i,1}}\,\Xi_{h_i}={\mathfrak m}_{Y,z_i}\Xi_{h_i}\supset{\rm Ker}(h_i:\Xi_{h_i}\to\Xi_{h_i}).$$
\par\nin Combined with (\hl{5.1.8}{5.1.8}), this implies for $\xi$ as in (\hl{5.1.6}{5.1.6})
\htt{5.2.1}{}
$$\bl[\tfrac{p}{d}\,\xi-\ddd\eta\br]\in V^{>\alpha_{h_i,1}\,}\Xi_{h_i}.
\leqno(5.2.1)$$
\par\nin So the infinite dimensionality of $M^{(2)}$ follows from (\hl{5.1.5}{5.1.5}), (\hl{5.1.7}{5.1.7}) and (\hl{5.2.1}{5.2.1}). It implies the assertion for $N^{(2)}$ since the morphisms
$$\ddd^{(1)}:N_{p+d}\to M_p$$
\par\nin are morphisms of finite dimensional vector spaces of the same dimension for $p\gg 0$. The assertion for the torsion part $(H^nA_f^{\ssb})_{\rm tor}$ then follows from Theorem~\hl{T2}{2}. This finishes the proof of Theorem~\hl{T5.2}{5.2}.
\par\htt{T5.3}{}\msn
{\bf Theorem~5.3.} {\it With the notation of {\rm\hl{5.1}{5.1}}, assume the $h_i$ are quasihomogeneous $($that is, $h_i\in(\dd h_i))$ for any $i\in I_0$. Then the kernel and cokernel of $d^{(1)}:N\to M$ $($that is, $N^{(2)}$ and $M^{(2)}$ in Theorem~$\hl{T2}{2})$ and $(H^nA_f^{\ssb})_{\rm tor}$ are finite dimensional over $\C$. More precisely we have
\htt{5.3.1}{}
$$\nu^{(2)}_{p+d}:=\dim N^{(2)}_{p+d}\les\#\bl\{(i,l)\,\big|\,\,\alpha_{h_i,l}=\tfrac{p}{d}\,\,(i\in I_0)\,\br\},
\leqno(5.3.1)$$
\par\nin and the equality holds in the case where $\mu'_p=0$ and either $\nu_{p+d}=\tauZ$ or all the singularities of $Z$ are ordinary double points.}
\msn
{\it Proof.} By Theorem~\hl{T2}{2}, it is enough to show the inequality (\hl{5.3.1}{5.3.1}) together with the equality in the special case as above.
Take any $i\in I_0$. In the notation of \hl{5.1}{5.1} there is a local analytic coordinate system $(y_1,\cdots,y_{n-1})$ of $Y'$ around $z_i$ together with positive rational numbers $w_1,\dots,w_{n-1}$ such that $h_i$ is a linear combination of monomials $\prod_jy_j^{m_j}$ with $\sum_jw_jm_j=1$ (see \cite{SaK}). Then
$$v(h_i)=h_i\q\h{with}\q v:=\msum_j\,w_j\,y_j\,\dd_{y_j}.$$
\par\nin We will denote the contraction of $\ddd y_1\wedge\dots\wedge y_{n-1}$ and $v$ by $\zeta$.
\sk
Take a monomial basis $\{\xi_l\}$ of $\Xi_{h_i}$, where monomial means that
$$\xi_l=\h{$\prod$}_j\,y_j^{m_{l,j}}\,\ddd y_1\wedge\dots\wedge\ddd y_{n-1}\q\h{with}\q m_{l,j}\in\N.$$
\par\nin Set
$$\eta_l:=\h{$\prod$}_j\,y_j^{m_{l,j}}\,\zeta,\q w(\xi_l):=\msum_j\,w_j(m_{l,j}+1).$$
\par\nin Then
$$\df\wedge\eta_l=f\,\xi_l,\q\ddd\eta_l=w(\xi_l)\,\xi_l.$$
\par\nin So we get
$$\dd_tt\,[\xi_l]=w(\xi_l)\,[\xi_l]\q\h{in}\,\,\,H''_{h_i}.$$
\par\nin In particular
$$w(\xi_l)=\alpha_{h_i,l},$$
\par\nin by changing the ordering of the $\xi_l$ if necessary.
The inequality (\hl{5.3.1}{5.3.1}) then follows from (\hl{5.1.7}{5.1.7}) and (\hl{5.1.8}{5.1.8}) together with the inclusion (\hl{5.1.2}{5.1.2}).
In case the assumption after (\hl{5.3.1}{5.3.1}) is satisfied, we have the equality by using the remark after (\hl{5.1.2}{5.1.2}) together with the assertion that $\alpha_{h_i,l}=(n{-}1)/2$ if $z_i$ is an ordinary double point of $Z$.
This finishes the proof of Theorem~\hl{T5.3}{5.3}.
\par\htt{C5.4}{}\msn
{\bf Corollary~5.4.} {\it With the hypothesis of Theorem~{\rm\hl{T5.3}{5.3}}, assume $n=2$ or more generally
\htt{5.4.1}{}
$$\max\bl\{\,\alpha_{h_i,l}\,\,\big|\,\,d\alpha_{h_i,l}\in\N\,\,\,(i\in I_0)\,\br\}<1+\tfrac{n}{d},
\leqno(5.4.1)$$
\par\nin $($for instance, $d\alpha_{h_i,l}\notin\N$ for any $l$ and $i\in I_0)$. Then the pole order spectral sequences $(\hl{4.4.4}{4.4.4})$ and $(\hl{4.5.1}{4.5.1})$ degenerate at $E_2$, and $(H^nA^{\ssb}_f)_{\rm tor}=0$.}
\msn
{\it Proof.} This follows from Theorem~\hl{T5.3}{5.3} together with Corollary~\hl{C4.7}{4.7} and Theorem~\hl{T2}{2} since $\ddd^{(r)}$ is a graded morphism of degree $-rd$.
\par\htt{C5.5}{}\msn
{\bf Corollary~5.5.} {\it With the first hypothesis of Theorem~{\rm\hl{T5.3}{5.3}}, assume $n\,{=}\,3$. Set $\alpha'_f\,{:=}\,\min\{\alpha_{h_i,l}\}$ in the notation of $(\hl{5.1.9}{5.1.9})$. Then}
\htt{5.5.1}{}
$$\nu_{p+d}=0\q\h{for}\q p<d\alpha'_f.
\leqno(5.5.1)$$
\par\nin \msn
{\it Proof.} Note first that $\alpha'_f\les 1$ since $\dim Z=1$. Assume $\nu_{p+d}\ne 0$ with $p<d\alpha'_f$. Then the image of $\ddd^{(1)}:N_{p+d}\to M_p$ is nonzero by Theorem~\hl{T5.3}{5.3}. Hence we get by Theorem~\hl{T2}{2}
$${}^P\!n_{f,p/d}<\mu_p\les\binom{p-1}{n-1},$$
\par\nin where the ${}^P\!n_{f,p/d}$ are the coefficients of the pole order spectrum $\Sp_P(f)$ as in \hl{3.2}{3.2}. However, this contradicts Remark~\hl{R3.6}{3.6}\,(ii). So Corollary~\hl{C5.5}{5.5} follows.
\par\htt{R5.6}{}\msn
{\bf Remarks~5.6.} (i) In Theorem~\hl{T5.3}{5.3}, the inequality (\hl{5.3.1}{5.3.1}) holds with the left-hand side replaced by the dimension of the kernel of the composition
$$N_{p+d}\buildrel{\ddd^{(1)}}\over\longrightarrow M_p\to M''_p.$$
\par\nin Indeed, (\hl{5.1.1}{5.1.1}) implies that (\hl{5.1.2}{5.1.2}) holds with $N_{p+d}^{(2)}$ replaced by this kernel.
\ms
(ii) Corollary~\hl{C5.4}{5.4} seems to be closely related to the short exact sequence in \cite[Thm.~1]{DiSa1}.
\ms
(iii) If all the singularities of $Z$ are nodes, then $\alpha'_f=1$, and the estimation obtained by Corollary~\hl{C5.5}{5.5} coincides with the one in \cite[Thm.~4.1]{DiSt2}, which is known to be sharp. It is also sharp if the singularities are $A_1$ or $D_4$ (for instance, $f=(x^2-y^2)(x^2-z^2)(y^2-z^2)$).
\ms
(iv) The proof of the finiteness of $(H^nA_f^{\ssb})_{\rm tor}$ can be reduced to the analytic local case by considering the formal completion where the direct sum is replaced with the infinite direct product and the convergent power series factors through the formal completion.
\ms
(v) If $n=3$ and $Z$ has only ordinary double points as singularities, then the coefficients $n_{f,\alpha}$ of the Steenbrink spectrum for $\alpha\notin\Z$ are the same as that of a central hyperplane arrangement in $\C^3$ having only ordinary double points in $\PP^2$. (Note that its formula can be found in \cite{BuSa}.)
Indeed, the vanishing cycle sheaf $\varphi_{f,\ne 1}\Q_X$ is supported at the origin so that we have the symmetry of the coefficients $n_{f,\alpha}$ for $\alpha\notin\Z$. Moreover $n_{f,\alpha}$ for $\alpha<1$ can be obtained by Remark~\hl{R3.6}{3.6}\,(ii), and $n_{f,\alpha}$ for $\alpha\in(1,2)$ can be calculated from the $n_{f,\alpha}$ for $\alpha\notin(1,2)$ by using the relation with the Euler characteristic of the complement of $Z\subset\PP^2$. (The latter follows from (\hl{3.1.2}{3.1.2}).)
Note also that the $n_{f,\alpha}$ for $\alpha\in\Z$ can be obtained from the Hodge numbers of the complement of $Z\subset\PP^{n-1}$.
\msn
{\bf Examples~5.7.}
We first give some examples where the assumptions of Corollary~\hl{C5.4}{5.4} and the last conditions of Theorem~\hl{T5.3}{5.3} are all satisfied, and moreover Remark~\hl{R5.6}{5.6}\,(v) can also be applied.
These are also examples of type (I) singularities (that is, (\hl{5}{5}) in the introduction is satisfied).
\skn
(i) $f=xyz\,\,\,$(three $A_1$ singularities in $\PP^2$) $\,\,n=d=3.$
$$\begin{array}{rccccccccccccccc}
k\,: &1 &2 &3 &4 &5 &6 &7 &8 &9 &\cdots \\
\gamma_k: & & &1 &3 &3 &1\\
\mu'_k: \\
\mu''_k: & & &1 &3 &3 &3 &3 &3 &3 &\cdots\\
\mu_k: & & &1 &3 &3 &3 &3 &3 &3 &\cdots\\
\nu_k: & & & & & &2 &3 &3 &3 &\cdots\\
\mu^{\scriptscriptstyle(2)}_k: & & &1\\
\nu^{\scriptscriptstyle(2)}_k: & & & & & &2\\
{}^P\!n_{f,k/d}: & & &1 & & &-2\\
n_{f,k/d}: & & &1 & & &-2\\
\end{array}$$
\par\nin \msn
(ii) $f=x^2y^2+x^2z^2+y^2z^2\,\,\,$(three $A_1$ singularities in $\PP^2$) $\,n=3,\,d=4.$
$$\begin{array}{rccccccccccccccccc}
k\,: &1 &2 &3 &4 &5 &6 &7 &8 &9 &10 &11 &12 &\cdots&\raise-3mm\h{ }\\
\gamma_k: & & &1 &3 &6 &7 &6 &3 &1\\
\mu'_k: & & & & &3 &4 &3\\
\mu''_k: & & &1 &3 &3 &3 &3 &3 &3 &3 &3 &3 &\cdots\\
\mu_k: & & &1 &3 &6 &7 &6 &3 &3 &3 &3 &3 &\cdots\\
\nu_k: & & & & & & & & &2 &3 &3 &3 &\cdots \\
\mu^{\scriptscriptstyle(2)}_k: & & &1 &3 &4 &4 &3\\
\nu^{\scriptscriptstyle(2)}_k: \\
{}^P\!n_{f,k/d}: & & &1 &3 &4 &4 &3 &0 &0\\
n_{f,k/d}: & & &1 &3 &3 &4 &3 &0 &1\\
\end{array}$$
\par\nin \msn
(iii) $f=xyz(x+y+z)\,\,\,$(six $A_1$ singularities in $\PP^2$) $\,n=3,\,d=4.$
$$\begin{array}{rccccccccccccccccc}
k\,: &1 &2 &3 &4 &5 &6 &7 &8 &9 &10 &11 &12 &\cdots&\raise-3mm\h{ }\\
\gamma_k: & & &1 &3 &6 &7 &6 &3 &1\\
\mu'_k: & & & & & &1 & & &\\
\mu''_k: & & &1 &3 &6 &6 &6 &6 &6 &6 &6 &6 &\cdots\\
\mu_k: & & &1 &3 &6 &7 &6 &6 &6 &6 &6 &6 &\cdots\\
\nu_k: & & & & & & & &3 &5 &6 &6 &6 &\cdots \\
\mu^{\scriptscriptstyle(2)}_k: & & &1 &3 &1 &1\\
\nu^{\scriptscriptstyle(2)}_k: & & & & & & & &3\\
{}^P\!n_{f,k/d}: & & &1 &3 &1 &1 &0 &-3 &0\\
n_{f,k/d}: & & &1 &3 &0 &1 &0 &-3 &1\\
\end{array}$$
\par\nin Here we have $\mu''_4=3$ by Lemma~\hl{L2.1}{2.1}, but the proof of $\mu''_5=6$ is not so trivial. Indeed, if $\mu''_5<6$, then we have $\nu_{\,7}\ne 0$ by Corollary~2. However, this contradicts Corollary~\hl{C5.5}{5.5}.
\msn
{\bf Examples~5.8.}
(i) $f=x^2y^2+z^4\,\,\,$(two $A_3$ singularities in $\PP^2$) $\,n=3,\,d=4.$
\msn
The calculation of this example does not immediately follow from Corollary~\hl{C5.4}{5.4}, since the last conditions of Theorem~\hl{T5.3}{5.3} are not satisfied and Remark~\hl{R5.6}{5.6}\,(v) does not apply to this example.
This example can be calculated by using the Thom-Sebastiani type theorems in Proposition~\hl{P2.2}{2.2} and Section~\hl{4.9}{4.9}.
\sk
$$\begin{array}{rccccccccccccccccc}
k\,: &1 &2 &3 &4 &5 &6 &7 &8 &9 &10 &11 &12 &\cdots&\raise-3mm\h{ }\\
\gamma_k: & & &1 &3 &6 &7 &6 &3 &1\\
\mu'_k: & & & & &1 &1 &1\\
\mu''_k: & & &1 &3 &5 &6 &6 &6 &6 &6 &6 &6 &\cdots\\
\mu_k: & & &1 &3 &6 &7 &7 &6 &6 &6 &6 &6 &\cdots\\
\nu_k: & & & & & & &1 &3 &5 &6 &6 &6 &\cdots \\
\mu^{\scriptscriptstyle(2)}_k: & & &1 &1 &2 &1 &1\\
\nu^{\scriptscriptstyle(2)}_k: & & & & & & &1 &1 &1\\\
{}^P\!n_{f,k/d}: & & &1 &1 &2 &1 &0 &-1 &-1\\
n_{f,k/d}: & & &1 &1 &2 &1 &0 &-1 &-1\\
\end{array}$$
\par\nin We note the calculation in the case $f=x^2y^2$ for the convenience of the reader.
\msn
(ii) $f=x^2y^2\,\,\,$(two $A_1$ singularities in $\PP^1$) $\,\,n=2,\,d=4.$
$$\begin{array}{rcccccccccccccccc}
k\,: &1 &2 &3 &4 &5 &6 &7 &8 &\cdots \\
\gamma_k: & &1 &2 &3 &2 &1\\
\mu'_k: & & & &1\\
\mu''_k: & &1 &2 &2 &2 &2 &2 &2 &\cdots\\
\mu_k: & &1 &2 &3 &2 &2 &2 &2 &\cdots\\
\nu_k: & & & & & &1 &2 &2 &\cdots\\
\mu^{\scriptscriptstyle(2)}_k: & &1 & &1\\
\nu^{\scriptscriptstyle(2)}_k: & & & & & &1\\
{}^P\!n_{f,k/d}: & &1 &0 &1 &0 &-1\\
n_{f,k/d}: & &1 &0 &1 &0 &-1\\
\end{array}$$
\par\nin \par\htt{R5.9}{}\msn
{\bf Remark~5.9.} We have the $V$-filtration of Kashiwara and Malgrange on $N_p$, $M''_p$ by using the injections in (5.1.1). Assume all the singularities of $Z$ are weighted homogeneous. It seems that the following holds in many examples:
\htt{5.9.1}{}
$$\dim\Gr_V^{\alpha}N_{p+d}=\begin{cases}\nu_{p+d}^{(2)}=\nu_{p+d}^{(\infty)}=n^1_{f,\alpha+1}&\h{if}\,\,\,p/d=\alpha,\\ 0&\h{if}\,\,\,p/d<\alpha,\end{cases}
\leqno(5.9.1)$$
\par\nin where $n^j_{f,\alpha}$ is as in (3.2.3), and $\nu_{p+d}^{(\infty)}:=\nu_{p+d}^{(r)}$ ($r\gg 0$).
Note that (\hl{5.9.1}{5.9.1}) would imply the $E_2$-degeneration of the pole order spectral sequence in Question~2, see also \cite{wh}.
\sk
As for $M''_p$, (\hl{5.9.1}{5.9.1}) seems to correspond by duality to the following:
\htt{5.9.2}{}
$$\dim\Gr_V^{\alpha}M''_p=\begin{cases}n_{Z,\alpha}-n^1_{f,\alpha}&\h{if}\,\,\,p/d=\alpha,\\ n_{Z,\alpha}&\h{if}\,\,\,p/d>\alpha,\end{cases}
\leqno(5.9.2)$$
\par\nin where $n_{Z,\alpha}:=\msum_{k\les r_0}\,n_{h_k,\alpha}$ with $n_{h_k,\alpha}$ defined for the isolated singularities $\{h_k=0\}$ ($k\les r_0$) as in (\hl{3.2.1}{3.2.1}).
Indeed, we have the symmetries
$$n_{Z,\alpha}=n_{Z,n-1-\alpha},\q n^1_{f,\alpha}=n^1_{f,n-\alpha}\q(\alpha\in\Q),$$
\par\nin and it is expected that the duality isomorphism in Theorem~\hl{T1}{1} is compatible with the filtration $V$ on $N_p$, $M''_p$ in an appropriate sense so that we have the equality
$$\dim\Gr_V^{\alpha}N_p+\dim\Gr_V^{n-1-\alpha}M''_{nd-p}=n_{Z,\alpha}\q(\alpha\in\Q,\,p\in\Z),
\leqno(5.9.3)$$
\par\nin giving a refinement of Corollary~\hl{C2}{2}. Note that (\hl{5.9.2}{5.9.2}) for $\alpha=p/d$ is closely related to \cite{Kl}.
\sk
If the above formulas hold, these would imply a refinement of Corollary~\hl{C5.5}{5.5} (and also its generalization to the case $n>3$ in \cite[Theorem~9]{DiSa3}).
However, it is quite nontrivial whether (\hl{5.9.2}{5.9.2}) holds, for instance, even for $p/d>\alpha$, since this is closely related to the independence of the $V$-filtrations associated to various singular points of $Z$.
In the case where the Newton boundaries of $f$ are non-degenerate, the formula for $M''_p$ with $\alpha\les 1$ seems to follow from the theories of multiplier ideals and microlocal $V$-filtrations.
\par\htt{5.10}{}\msn
{\bf 5.10.~Proof of Theorem~\hl{T3}{3}.} Set $\beta:=\max\Rc_Z$. This coincides with the maximum of the spectral numbers of the singularities of $Z$, since the singularities are weighted homogeneous. By Theorem~\hl{T5.3}{5.3} we then get that
\htt{5.10.1}{}
$$N^{(2)}_{k+d}=0\q\h{if}\q\tfrac{k}{d}>\beta.
\leqno(5.10.1)$$
\par\nin This implies the partial $E_2$-degeneration
\htt{5.10.2}{}
$$M^{(2)}_k=M^{(r)}_k\q\h{if}\q\tfrac{k}{d}>\beta{-}1,\,r>2,
\leqno(5.10.2)$$
\par\nin since the differential $\ddd^{(r)}:N^{(r)}_k\to M^{(r)}_{k-rd}$ shifts the degree by $-rd$, see Theorem~\hl{T2}{2}.
\sk
On the other hand, we get by the using symmetry in Corollary~\hl{C2}{2}
\htt{5.10.3}{}
$$\delta''_k=0\q\h{if}\q k>(n{-}1)d{-}n.
\leqno(5.10.3)$$
\par\nin Indeed, if $\nu_{d+n-1}\ne0$, then $f$ is annihilated by a linear combination of the $\dd_{x_i}$, and is a polynomial of $n{-}1$ variables after a linear coordinate change.
\sk
Now let $k$ be as in Theorem~\hl{T3}{3}. The above arguments imply that
\htt{5.10.4}{}
$$\dim M^{(\infty)}_k=\dim M'_k.
\leqno(5.10.4)$$
\par\nin Indeed, $\beta<n{-}1$ by \cite{mic}, and the injectivity of $d^{(1)}_k:N_{k+d}\to M_k$ follows from Theorem~\hl{T5.3}{5.3} using the hypothesis of Theorem~\hl{T3}{3} or \cite{dFEM} in the case $n=3$. The assertion then follows from \cite[Theorem 2]{cm-b} (since $\beta<n{-}1$). This finishes the proof of Theorem~\hl{T3}{3}.

\end{document}